\input amstex
\documentstyle{amsppt}
%
\catcode`@=11
\redefine\output@{%
  \def\break{\penalty-\@M}\let\par\endgraf
  \ifodd\pageno\global\hoffset=105pt\else\global\hoffset=8pt\fi  
  \shipout\vbox{%
    \ifplain@
      \let\makeheadline\relax \let\makefootline\relax
    \else
      \iffirstpage@ \global\firstpage@false
        \let\rightheadline\frheadline
        \let\leftheadline\flheadline
      \else
        \ifrunheads@ 
        \else \let\makeheadline\relax
        \fi
      \fi
    \fi
    \makeheadline \pagebody \makefootline}%
  \advancepageno \ifnum\outputpenalty>-\@MM\else\dosupereject\fi
}
\def\Beta{\mathchar"0\hexnumber@\rmfam 42}
\catcode`\@=\active
\nopagenumbers
\chardef\textvolna='176

\chardef\bigalpha='013
\def\negskp{\hskip -2pt}

\def\Img{\operatorname{Im}}

\chardef\degree="5E
\def\const{\operatorname{const}}
\def\blue#1{#1}

\gdef\darkred#1{#1}
\catcode`#=11\def\diez{#}\catcode`#=6
\catcode`&=11\catcode`&=4
\catcode`_=11\def\podcherkivanie{_}\catcode`_=8
\catcode`\^=11\catcode`\^=7
\catcode`~=11\def\volna{~}\catcode`~=\active
\def\mycite#1{\cite{\blue{#1}}\immediate\special{ps:
     ShrHPSdict begin /ShrBORDERthickness 0 def}}
\def\myciterange#1#2#3#4{\cite{\blue{#2#3#4}}\immediate\special{ps:
     ShrHPSdict begin /ShrBORDERthickness 0 def}}
\def\mytag#1{%
    \tag#1}
\def\mythetag#1{\thetag{\blue{#1}}\immediate\special{ps:
     ShrHPSdict begin /ShrBORDERthickness 0 def}}
\def\myrefno#1{\no#1}
\def\myhref#1#2{\blue{#2}\immediate\special{ps:
     ShrHPSdict begin /ShrBORDERthickness 0 def}}
\def\myEarXivlink{\myhref{http://arXiv.org}{http:/\negskp/arXiv.org}}

\def\mytheorem#1{\csname proclaim\endcsname{Theorem #1}}
\def\mytheoremwithtitle#1#2{\csname proclaim\endcsname{Theorem #1#2}}
\def\mythetheorem#1{\blue{#1}\immediate\special{ps:
     ShrHPSdict begin /ShrBORDERthickness 0 def}}
\def\mylemma#1{\csname proclaim\endcsname{Lemma #1}}
\def\mylemmawithtitle#1#2{\csname proclaim\endcsname{Lemma #1#2}}

\def\mycorollary#1{\csname proclaim\endcsname{Corollary #1}}

\def\mydefinition#1{\definition{Definition #1}}
\def\mythedefinition#1{\blue{#1}\immediate\special{ps:
     ShrHPSdict begin /ShrBORDERthickness 0 def}}
\def\myconjecture#1{\csname proclaim\endcsname{Conjecture #1}}
\def\myconjecturewithtitle#1#2{\csname proclaim\endcsname{Conjecture #1#2}}

\def\myproblem#1{\csname proclaim\endcsname{Problem #1}}
\def\myproblemwithtitle#1#2{\csname proclaim\endcsname{Problem #1#2}}


\pagewidth{360pt}
\pageheight{606pt}
\topmatter
\title
Asymptotic estimates for roots of the cuboid characteristic 
equation in the nonlinear region.
\endtitle
\rightheadtext{Asymptotic estimates for roots \dots}
\author
Ruslan Sharipov
\endauthor
\address Bashkir State University, 32 Zaki Validi street, 450074 Ufa, Russia
\endaddress
\email
\myhref{mailto:r-sharipov\@mail.ru}{r-sharipov\@mail.ru}
\endemail
\abstract
     A perfect cuboid is a rectangular parallelepiped. Its edges, its face 
diagonals, and its space diagonal are of integer lengths. None of such cuboids
is known thus far, though the system of Diophantine equations describing them 
is easily written. The cuboid characteristic equation is a twelfth degree 
Diophantine equation derived from the initial cuboid equations and equivalent 
to them. In the case of the second cuboid conjecture it reduces to a tenth 
degree equation. This equation comprises two parameters. Previously various 
asymptotics for roots of this equation were studied as its parameters tend 
to infinity either separately or simultaneously provided some linear combination
of them is preserved finite. In the present paper this linear combination is
replaced by a certain nonlinear expression. 
\endabstract
\subjclassyear{2000}
\subjclass 11D41, 11D72, 41A60\endsubjclass
\endtopmatter
\TagsOnRight
\document

\head
1. Introduction.
\endhead
     Omitting details, which can be found in \myciterange{1}{1}{--}{3}, let's 
proceed to the tenth degree cuboid characteristic equation arising in the case 
of the second cuboid conjecture:
$$
\hskip -2em
Q_{p\kern 0.6pt q}(t)=0.
\mytag{1.1}
$$
The tenth degree polynomial $Q_{p\kern 0.6pt q}(t)$ in \mythetag{1.1} is given
by the explicit formula 
$$
\gathered
Q_{p\kern 0.6pt q}(t)=t^{10}+(2\,q^{\kern 0.7pt 2}+p^{\kern 1pt 2})\,(3
\,q^{\kern 0.7pt 2}-2\,p^{\kern 1pt 2})\,t^8+(q^{\kern 0.5pt 8}+10
\,p^{\kern 1pt 2}\,q^{\kern 0.5pt 6}+\\
+\,4\,p^{\kern 1pt 4}\,q^{\kern 0.5pt 4}-14\,p^{\kern 1pt 6}\,q^{\kern 0.7pt 2}
+p^{\kern 1pt 8})\,t^6-p^{\kern 1pt 2}\,q^{\kern 0.7pt 2}
\,(q^{\kern 0.5pt 8}-14\,p^{\kern 1pt 2}\,q^{\kern 0.5pt 6}
+4\,p^{\kern 1pt 4}\,q^{\kern 0.5pt 4}+\\
+\,10\,p^{\kern 1pt 6}\,q^{\kern 0.7pt 2}+p^{\kern 1pt 8})\,t^4
-p^{\kern 1pt 6}\,q^{\kern 0.5pt 6}\,(q^{\kern 0.7pt 2}
+2\,p^{\kern 1pt 2})\,(3\,p^{\kern 1pt 2}-2\,q^{\kern 0.7pt 2})\,t^2
-q^{\kern 0.7pt 10}\,p^{\kern 1pt 10}.
\endgathered\quad
\mytag{1.2}
$$
More details concerning the polynomial \mythetag{1.2} and its background can
be found in \myciterange{4}{4}{--}{8}. For the history and various 
approaches to the problem of perfect cuboids the reader is referred to \myciterange{9}{9}{--}{55}. The papers \myciterange{56}{56}{--}{68} constitute
a separate stream using a symmetry approach to the perfect cuboid equations. 
This approach is different from the approach of the present paper. Therefore 
we do not consider the results of the papers \myciterange{56}{56}{--}{68} 
below.\par
    The tenth degree polynomial \mythetag{1.2} is related to the perfect cuboid 
problem through the following theorem (see Theorem~8.1 in \mycite{1} or in
\mycite{2}).\par
\mytheorem{1.1} A triple of positive integer numbers $p$, $q$, and $t$ satisfying 
the equation \mythetag{1.1} and such that $p\neq q$ are coprime produces a perfect
cuboid if and only if the following inequalities are
fulfilled: 
$$
\xalignat 4
&t>p^{\kern 1pt 2},
&&t>p\,q,
&&t>q^{\kern 0.7pt 2},
&&(p^{\kern 1pt 2}+t)\,(p\,q+t)>2\,t^2.
\endxalignat
$$
\endproclaim
     The mechanism associating the numbers $p,\,q,\,t$ with perfect cuboids was
found in \mycite{4} and \mycite{5}. It is described in brief in \mycite{1},
\mycite{2}, and \mycite{3}. We shall not reproduce this 
\vadjust{\vskip 406pt\hbox to 0pt{\kern 10pt
\includegraphics{Strategy04.eps}\hss}\vskip 0pt}description in the 
present paper.\par
     Relying on Theorem~\mythetheorem{1.1}, in \mycite{1}, \mycite{2} two forms 
of asymptotics for the roots of the tenth degree polynomial equation \mythetag{1.1} 
were studied:
$$
\pagebreak 
\xalignat 2
&\hskip -2em
p=\const,&&q\to+\infty,\\
\vspace{-1.5ex}
\mytag{1.3}\\
\vspace{-1.5ex}
&\hskip -2em
q=\const,&&p\to+\infty.
\endxalignat
$$ 
As a result of studying limits \mythetag{1.3} in paper \mycite{2} the 
following three regions in the positive quadrant of the $p\,q$\,-\,coordinate 
plane were defined:
\roster
\item"1)" {\bf linear region} given by the linear inequalities
$$
\xalignat 2
&\frac{q}{59}<p,
&&p<59\,q;
\mytag{1.4}
\endxalignat
$$
\item"2)" {\bf nonlinear region} given by the nonlinear inequalities
$$
\xalignat 2
&59\,q\leqslant p,
&&p\leqslant 9\,q^{\kern 0.7pt 3};
\mytag{1.5}
\endxalignat
$$
\item"3)" {\bf no cuboid region} which is the rest of the positive 
$p\,q$\,-\,quadrant.
\endroster
These regions are schematically shown in Fig\.~1.1 above. The linear region 
\mythetag{1.4} is shown in sky blue. The nonlinear region \mythetag{1.5} is shown 
in faded pink. And finally, the no cuboid region is shown in gray.\par
     The linear region was studied in \mycite{3}. As a result a narrow
strip surrounding the bisector line $p=q$ was cut off from the linear region. 
It is given by the inequalities 
$$
\xalignat 2
&\hskip -2em
q-\frac{q}{97}\leqslant p,
&&p\leqslant q+\min\Bigl(\frac{q}{97},\root3\of{\frac{q}{74}}\Bigr).
\mytag{1.6}
\endxalignat 
$$
The strip \mythetag{1.6} was annexed to the no cuboid region. Unfortunately
the rest of the linear region \mythetag{1.4} still remains uncertain. Perfect
cuboids can potentially be found in it, but none of them is actually found.
The main goal of the present paper is to study the nonlinear region 
\mythetag{1.5}.\par  
\head
2. Nonlinear transformation of parameters.
\endhead
     The upper boundary of the nonlinear region is given by the cubic parabola
$p=9\,q^{\kern 0.7pt 3}$ (see \mythetag{1.5} and Fig\.~1.1 above). For this reason
we consider the following cubic transformation of the parameters $p$ and $q$:
$$
\xalignat 2
&\hskip -2em
	\tilde p=B\,q^{\kern 0.7pt 3}-p,
&&\tilde q=q.
\mytag{2.1}
\endxalignat
$$
Here $B$ is some positive integer number. The transformation \mythetag{2.1} 
is invertible:
$$
\xalignat 2
&\hskip -2em
p=B\,\tilde q^{\kern 0.7pt 3}-\tilde p,
&&q=\tilde q.
\mytag{2.2}
\endxalignat
$$
The transformations \mythetag{2.1} and the transformation \mythetag{2.2} map
the integer $p\,q$\,-\,grid onto the integer $\tilde p\,\tilde q$\,-\,grid
and vice versa.\par
     Note that the curve given by the condition $\tilde p=\const$ is a cubic 
parabola. In particular, if $B=9$ and $\tilde p=0$, it coincides with the upper
boundary of the nonlinear region in Fig\.~1.1. For this reason, instead of 
\mythetag{1.3}, we set
$$
\xalignat 2
&\hskip -2em
\tilde p=\const,&&\tilde q\to+\infty
\mytag{2.3}
\endxalignat
$$\par
     In order to study the limit \mythetag{2.3} let's substitute \mythetag{2.2}
into the polynomial \mythetag{1.2}. As a result we get another polynomial 
$Q_{\tilde p\kern 0.6pt \tilde q}(t)$. This polynomial is given by an explicit 
formula. However, the formula for the polynomial 
$Q_{\tilde p\kern 0.6pt \tilde q}(t)$ is rather huge. In the fully expanded form 
it has $108$ summands. For this reason it is placed to the ancillary file 
\darkred{{\tt strategy\kern -0.5pt\_\kern 1.5pt formulas\_\kern 0.5pt 04.txt}} 
in a machine-readable form.\par
     Using the polynomial $Q_{\tilde p\kern 0.6pt \tilde q}(t)$, one can write
an equation similar to \mythetag{1.1}:
$$
\hskip -2em
Q_{\tilde p\kern 0.6pt \tilde q}(t)=0.
\mytag{2.4}
$$
It is clear that the equation \mythetag{2.4} has the same roots as the 
equation \mythetag{1.1}, though they are expressed through different parameters. 
Like \ $Q_{p\kern 0.6pt q}(t)$, the polynomial $Q_{\tilde p\kern 0.6pt \tilde q}(t)$
in \mythetag{2.4} is even with respect to $t$. Along with each root $t$ it has 
the opposite root $-t$. Therefore we use the condition
$$
\hskip -2em
\cases \text{$t>0$ \ if \ $t$ \ is a real root,}\\
\text{$\Img(t)>0$ \ if \ $t$ \ is a complex root,}
\endcases
\mytag{2.5}
$$
in order to divide the roots of the equation \mythetag{2.4} into two groups. 
We denote through $t_1$, $t_2$, $t_3$, $t_4$, $t_5$ the roots that obey the 
conditions \mythetag{2.5}. Then $t_6$, $t_7$, $t_8$, $t_9$, $t_{10}$ are 
opposite roots of the equation \mythetag{2.4}:
$$
\xalignat 5
&t_6=-t_1,&&t_7=-t_2,&&t_8=-t_3,&&t_9=-t_4,&&t_{10}=-t_5.
\qquad\quad
\endxalignat
$$\par
\head
3. Parabolic expansions. 
\endhead
    Typically, asymptotic expansions for roots of a polynomial equation look 
like power series (see \mycite{69}). By analogy to (2.3) in \mycite{1} and
according to \mythetag{2.3}, we write
$$
\hskip -2em
t_i(\tilde p,\tilde q)=C_i\,\tilde q^{\,\alpha_i}
\biggl(1+\sum^\infty_{s=1}\beta_{is}\,\tilde q^{-s}\biggr)
\text{\ \ as \ }\tilde q\to+\infty.
\mytag{3.1}
$$
The coefficients $C_i$ in \mythetag{3.1} should be nonzero: 
$C_i\neq 0$. The expansions \mythetag{3.1} here are called parabolic 
expansions since they occur along cubic parabolas in the plane of the 
original parameters $p$ and $q$.\par
     The main tool for studying asymptotic expansions of the form 
\mythetag{3.1} for roots of a polynomial equation is the Newton polygon
of the corresponding polynomial.  
\mydefinition{3.1} For any polynomial of two variables $P(t,q)$ the convex 
hull of all integer nodes $(m,r)$ on the coordinate plane associated with 
nonzero monomials $A_{m\kern 1pt r}\,q^{\kern 1pt r}\,t^m$ of this polynomial 
is called the Newton polygon of $P(t,q)$.  
\enddefinition
{\bf Remark}. Note that in our case the polynomial $Q_{\tilde p\kern 0.6pt 
\tilde q}(t)$ depends on three variables $\tilde p$, $\tilde q$, and 
$\tilde t$. However, we treat $\tilde p$ as a parameter and consider 
$Q_{\tilde p\kern 0.6pt \tilde q}(t)$ as a polynomial of two variables 
when applying Definition~\mythedefinition{3.1} to it.\par
     The Newton polygon of $Q_{\tilde p\kern 0.6pt \tilde q}(t)$ is shown 
in Fig\.~3.1 below. Its boundary consists of three parts --- the upper 
part, the lower part, and the vertical part. The upper part is drawn in green, 
the lower part is drawn in red. Here are the coefficients of monomials associated 
with the nodes on the upper part of the boundary:
$$
\xalignat 3
&\hskip -2em
A_{\kern 0.5pt 0\kern 2pt 40}=-B^{\kern 1pt 10},
&&A_{\kern 0.5pt 2\kern 1.5pt 36}=-6\,B^{\kern 1pt 10},
&&A_{\kern 0.5pt 4\kern 1.5pt 32}=-B^{\kern 1pt 10},
\quad
\\
\vspace{-1.5ex}
\mytag{3.2}\\
\vspace{-1.5ex}
&\hskip -2em
A_{\kern 0.5pt 6\kern 2pt 24}=B^{\kern 1pt 8},
&&A_{\kern 0.5pt 8\kern 1.5pt 12}=-2\,B^{\kern 1pt 4},
&&A_{\kern 0.5pt 10\kern 2pt 0}=1.
\quad
\endxalignat
$$
\mytheorem{3.1} The values of exponents $\alpha_i$ in the expansion \pagebreak
\mythetag{3.1} for roots of the equation \mythetag{2.4} are determined according 
to the formula $\alpha_i=-k$, where $k$ stands for slopes of segments of the 
polygonal line being the upper boundary of the Newton polygon in Fig\.~3.1. 
\endproclaim
      Theorem~\mythetheorem{3.1} is a standard fact in Newton polygons 
application to studying bivariate polynomials. \vadjust{\vskip 366pt\hbox 
to 0pt{\kern 10pt\includegraphics{Strategy05.eps}\hss}\vskip 0pt}Its 
proof was given in \mycite{1} for the sake of reader's convenience. In our 
particular case, applying Theorem~\mythetheorem{3.1}, we get 
$$
\xalignat 3
&\hskip -2em
\alpha_i=2,
&&\alpha_i=4,
&&\alpha_i=6.
\mytag{3.3}
\endxalignat
$$
The options \mythetag{3.3} determine the growth rates for the roots of the
equation \mythetag{2.4} as $\tilde q\to+\infty$. They grow as the second, 
the fourth, and the sixth powers of $\tilde q$.\par 
    {\bf The case }$\alpha_i=2$. This case corresponds to the topmost
segment on the upper boundary of the Newton polygon in Fig\.~3.1. This 
segment comprises three nodes $A_{\kern 0.5pt 0\kern 1.5pt 40}$,
$A_{\kern 0.5pt 2\kern 1.5pt 36}$, and $A_{\kern 0.5pt 4\kern 2pt 32}$. 
Therefore, substituting the expansion \mythetag{3.1} with $\alpha_i=2$ 
into the equation \mythetag{2.4}, we get the following equation for $C_i$:
$$
\hskip -2em
A_{\kern 0.5pt 4\kern 1.5pt 32}\ {C_i}^4
+A_{\kern 0.5pt 2\kern 1.5pt 36}\ {C_i}^2
+A_{\kern 0.5pt 0\kern 2pt 40}=0. 
\mytag{3.4}
$$
Taking into account \mythetag{3.2}, the equation \mythetag{3.4} is 
transformed to
$$
\hskip -2em
B^{\kern 1pt 10}\,{C_i}^4+6\,B^{\kern 1pt 10}\,{C_i}^2+B^{\kern 1pt 10}=0.
\mytag{3.5}
$$
Since $B\neq 0$ in \mythetag{2.1} and \mythetag{2.2}, we can cancel 
$B^{\kern 1pt 10}$ in \mythetag{3.5} and write \mythetag{3.5} as 
$$
\hskip -2em
{C_i}^4+6\,{C_i}^2+1=0.
\mytag{3.6}
$$
The equation \mythetag{3.6} has four purely complex roots:
$$
\xalignat 2
&\hskip -2em
C_i=(\sqrt{2}+1)\,\goth i,
&&C_i=(\sqrt{2}-1)\,\goth i,
\mytag{3.7}\\
&\hskip -2em
C_i=-(\sqrt{2}+1)\,\goth i,
&&C_i=-(\sqrt{2}-1)\,\goth i.
\mytag{3.8}
\endxalignat
$$
Here $\goth i=\sqrt{-1}$. The roots \mythetag{3.8} are excluded by the
condition \mythetag{2.5}. The remain is two root \mythetag{3.7} of 
multiplicity $1$. They yield the following asymptotic expansions:
$$
\hskip -2em
\aligned
&t_i(\tilde p,\tilde q)=(\sqrt{2}+1)\,\goth i\,\tilde q^{\kern 0.7pt 2}
\biggl(1+\sum^\infty_{s=1}\beta_{is}\,\tilde q^{\,-s}\biggr),\\
&t_i(\tilde p,\tilde q)=(\sqrt{2}-1)\,\goth i\,\tilde q^{\kern 0.7pt 2}
\biggl(1+\sum^\infty_{s=1}\beta_{is}\,\tilde q^{\,-s}\biggr). 
\endaligned
\mytag{3.9}
$$\par
     {\bf The case} $\alpha_i=4$. This case corresponds to the middle segment 
in the upper boundary of the Newton polygon in Fig\.~3.1. It comprises
two nodes $A_{\kern 0.5pt 4\kern 1.5pt 32}$ and $A_{\kern 0.5pt 6\kern 2pt 24}$. 
Therefore, substituting the expansion \mythetag{3.1} with $\alpha_i=4$ 
into the equation \mythetag{2.4}, we get the following equation for $C_i$:
$$
\hskip -2em
A_{\kern 0.5pt 6\kern 2pt 24}\ {C_i}^6
+A_{\kern 0.5pt 4\kern 1.5pt 32}\ {C_i}^4=0.
\mytag{3.10}
$$
Applying \mythetag{3.2}, the equation \mythetag{3.10} is transformed to
$$
\hskip -2em
B^{\kern 1pt 8}\,{C_i}^6-B^{\kern 1pt 10}\,{C_i}^4=0.
\mytag{3.11}
$$
Since $C_i\neq 0$ in \mythetag{3.1} and since $B\neq 0$ in \mythetag{2.1} and 
\mythetag{2.2}, we can cancel the common divisor ${B^{\kern 1pt 8}\,C_i}^4$ of two 
terms in \mythetag{3.11}. As a result \mythetag{3.11} takes the form 
$$
\hskip -2em
{C_i}^2-B^{\kern 1pt 2}=0.
\mytag{3.12}
$$
The equation \mythetag{3.12} has two real roots, which are simple:
$$
\xalignat 2
&\hskip -2em
C_i=B,
&&C_i=-B.
\mytag{3.13}
\endxalignat
$$
The second root \mythetag{3.13} is excluded by the condition \mythetag{2.5}. 
The remain is one simple positive root. It yields the following asymptotic 
expansion:
$$
\hskip -2em
t_i(\tilde p,\tilde q)=B\,\tilde q^{\kern 0.7pt 4}
\biggl(1+\sum^\infty_{s=1}\beta_{is}\,\tilde q^{\,-s}\biggr).
\mytag{3.14}
$$\par
     {\bf The case} $\alpha_i=6$.  This case corresponds to the lowermost 
segment in the upper boundary of the Newton polygon in Fig\.~3.1. It 
comprises three nodes $A_{\kern 0.5pt 6\kern 1.5pt 24}$,
$A_{\kern 0.5pt 8\kern 2pt 12}$,  and  $A_{\kern 0.5pt 10\kern 2pt 0}$.
Therefore, substituting the expansion \mythetag{3.1} with $\alpha_i=6$ 
into the equation \mythetag{2.4}, we get the following equation for $C_i$:
$$
\hskip -2em
A_{\kern 0.5pt 10\kern 2pt 0}\ {C_i}^{10}
+A_{\kern 0.5pt 8\kern 1.5pt 12}\ {C_i}^8
+A_{\kern 0.5pt 6\kern 1.5pt 24}\ {C_i}^6=0.
\mytag{3.15}
$$
Applying the formulas \mythetag{3.2}, the equation \mythetag{3.15} is 
transformed to
$$
\hskip -2em
{C_i}^{10}-2\,B^{\kern 1pt 4}\,{C_i}^8
+B^{\kern 1pt 8}\,{C_i}^6=0.
\mytag{3.16}
$$
Since $C_i\neq 0$ in \mythetag{3.1}, we can cancel the common divisor 
${C_i}^6$ of the three terms in \mythetag{3.16}. As a result the equation 
\mythetag{3.16} takes the form 
$$
\hskip -2em
{C_i}^4-2\,B^{\kern 1pt 4}\,{C_i}^2+B^{\kern 1pt 8}=0.
\mytag{3.17}
$$
The equation \mythetag{3.17} has two real roots of multiplicity $2$:
$$
\xalignat 2
&\hskip -2em
C_i=B^{\kern 1pt 2},
&&C_i=-B^{\kern 1pt 2}.
\mytag{3.18}
\endxalignat
$$
The second root \mythetag{3.18} is excluded by the condition \mythetag{2.5}. 
The remain is one simple positive root. It yields the following asymptotic 
expansion:
$$
\hskip -2em
t_i(\tilde p,\tilde q)=B^{\kern 1pt 2}\,\tilde q^{\kern 0.7pt 6}
\biggl(1+\sum^\infty_{s=1}\beta_{is}\,\tilde q^{\,-s}\biggr).
\mytag{3.19}
$$\par
The results \mythetag{3.9}, \mythetag{3.14}, \mythetag{3.19} are summed up
in the following theorem. 
\mytheorem{3.2} For sufficiently large positive values of the parameter 
$\tilde q$, i\.\,e\. for $\tilde q>\tilde q_{\text{\,min}}$, the tenth-degree 
equation \mythetag{2.4} has five roots of multiplicity one satisfying the 
condition \mythetag{2.5}. Three of them $t_1$, $t_2$, and $t_3$ are real roots. 
Their asymptotics as $\tilde q\to +\infty$ are given by the formulas 
$$
\xalignat 3
&\hskip -2em
t_1\sim B^{\kern 1pt 2}\,\tilde q^{\kern 0.7pt 6},
&&t_2\sim B^{\kern 1pt 2}\,\tilde q^{\kern 0.7pt 6},
&&t_3\sim B\,\tilde q^{\kern 0.7pt 4}.
\mytag{3.20}
\endxalignat
$$
The rest two roots $t_4$ and $t_5$ of the equation \mythetag{2.4} are complex. 
Their asymptotics as $\tilde q\to +\infty$ are given by the formulas 
$$
\xalignat 2
&\hskip -2em
t_4\sim (\sqrt{2}+1)\,\goth i\,\tilde q^{\kern 0.7pt 2},
&&t_5\sim (\sqrt{2}-1)\,\goth i\,\tilde q^{\kern 0.7pt 2}.
\mytag{3.21}
\endxalignat
$$
\endproclaim
The complex roots \mythetag{3.21} do not provide perfect cuboids. However, they 
are important for determining the exact number of real roots in asymptotic
intervals.\par
\head
4. Asymptotic estimates for real roots.
\endhead
     Acting by analogy to \mycite{1}, we replace the series expansion of the form 
\mythetag{3.1} by finite sum expansions with remainder terms. In the case of the
fast growing root $t_1$ we replace the expansion \mythetag{3.19} by the following 
sum:  
$$
\hskip -2em
\aligned
t_1=B^{\kern 1pt 2}\,\tilde q^{\kern 0.7pt 6}+2\,B\,\tilde q^{\kern 0.7pt 4}
&-2\,B\,\tilde p\,\tilde q^{\kern 1pt 3}-2\,\tilde q^{\kern 0.7pt 2}\,-\\
&-\,2\,\tilde p\,\tilde q+\tilde p^{\kern 1pt 2}+\frac{5}{B}
-\frac{20}{B^{\kern 1pt 2}\,\tilde q^{\kern 0.7pt 2}}
+R_1(\tilde p,\tilde q,B).
\endaligned
\mytag{4.1}
$$
The formula \mythetag{4.1} is in agreement with \mythetag{3.19}. 
It means that
$$
\xalignat 4
&\beta_{11}=0,
&&\beta_{12}=\frac{2}{B},
&&\beta_{13}=-\frac{2\,\tilde p}{B},
&&\beta_{14}=-\frac{2}{B^{\kern 1pt 2}},\\
&\beta_{15}=-\frac{2\,\tilde p}{B^{\kern 1pt 2}},
&&\beta_{16}=\frac{\tilde p^{\kern 1pt 2}}{B^{\kern 1pt 2}}
+\frac{5}{B^{\kern 1pt 3}},
&&\beta_{17}=0
&&\beta_{18}=-\frac{20}{B^{\kern 1pt 4}}
\endxalignat
$$
in the formula \mythetag{3.19}. Our goal is to derive an estimate 
of the form 
$$
\hskip -2em
|R_1(\tilde p,\tilde q,B)|<\frac{C}{\tilde q^{\kern 0.7pt 3}}
\mytag{4.2}
$$
for the remainder term in \mythetag{4.1}. In order to get such an estimate 
we substitute 
$$
t=B^{\kern 1pt 2}\,\tilde q^{\kern 0.7pt 6}
+2\,B\,\tilde q^{\kern 0.7pt 4}-2\,B\,\tilde p\,\tilde q^{\kern 1pt 3}
-2\,\tilde q^{\kern 0.7pt 2}-2\,\tilde p\,\tilde q
+\tilde p^{\kern 1pt 2}+\frac{5}{B}
-\frac{20}{B^{\kern 1pt 2}\,\tilde q^{\kern 0.7pt 2}}
+\frac{c}{\tilde q^{\kern 0.7pt 3}}.
\mytag{4.3}
$$
into the equation \mythetag{2.4}. Immediately after that we perform another
substitution into the equation obtained by substituting \mythetag{4.3}
into \mythetag{2.4}: 
$$
\hskip -2em
\tilde q=\frac{1}{z}.
\mytag{4.4}
$$
Upon two substitutions \mythetag{4.3} and \mythetag{4.4} and upon removing 
denominators the equation \mythetag{2.4} is written as a polynomial equation 
in the new variables $c$ and $z$:
$$
\hskip -2em
16\,c\,B^{\kern 1pt 37}=80\,\tilde p\,B^{\kern 1pt 35}+\varphi(c,z,\tilde p,B).
\mytag{4.5}
$$
Here $\varphi(c,z,\tilde p,B)$ is a polynomial of $c$, $z$, $\tilde p$, and 
$B$ with integer coefficients. The explicit expression for $\varphi(c,z,\tilde p,B)$
comprises $1612$ monomials. Therefore it is placed to the ancillary file 
\darkred{{\tt strategy\kern -0.5pt\_\kern 1.5pt formulas\_\kern 0.5pt 04.txt}} 
in a machine-readable form.\par
     Let's recall that $B$ is a positive integer constant. If $B\geqslant 10$,
then the curve $\tilde p=\const$ goes to infinity outside the nonlinear region 
(see \mythetag{2.1}, \mythetag{2.3}, and Fig\.~1.1). Therefore $B$ takes the 
following finite set of values:
$$
\hskip -2em
B=1,\,2,\,3,\,\ldots,\,9. 
\mytag{4.6}
$$
For each $B$ in \mythetag{4.6} assume that $c$ obeys the condition
$$
\hskip -2em
\aligned
-\frac{10\,|\tilde p\kern 0.7pt|}{B^{\kern 1pt 2}}<c<0\text{\ \ if \ }
\tilde p<0,\\
0<c<\frac{10\,|\tilde p\kern 0.7pt|}{B^{\kern 1pt 2}}\text{\ \ if \ }
\tilde p>0.
\endaligned
\mytag{4.7}
$$ 
The case $\tilde p=0$ is exceptional. It should be studied separately.\par
     Since $\tilde p$ is integer and since $\tilde p\neq 0$ in our present case,
we have the inequality 
$$
\hskip -2em
|\tilde p\kern 0.7pt|\geqslant 1.
\mytag{4.8}
$$
Since we study the asymptotics of the roots $t_i$ as $\tilde q\to+\infty$, we 
assume that 
$$
\hskip -2em
\tilde q\geqslant 20\,\root{3}\of{|\tilde p\kern 0.7pt|}.
\mytag{4.9}
$$
Let's apply \mythetag{4.9} to \mythetag{4.4}. As  a result we
derive the inequality $|z|\leqslant 1/20\,|\tilde p\kern 0.7pt|^{-1/3}$. 
Applying this inequality along with the inequalities \mythetag{4.7} and 
\mythetag{4.8} to the polynomial $\varphi(c,z,\tilde p,B)$, we derive the 
following estimate for it: 
$$
\hskip -2em
|\varphi(c,z,\tilde p,B)|\leqslant 72\,|\tilde p\kern 0.7pt|
\,B^{\kern 1pt 35}.
\mytag{4.10}
$$
The inequality \mythetag{4.10} holds for each value of $B$ in \mythetag{4.6}.
\par
     Now we apply the inequality \mythetag{4.10} to the equation \mythetag{4.5}.
If $\tilde p<0$, it means that the right hand side of the equation \mythetag{4.5} 
is a continuous function of $c$ that varies from $-152\,|\tilde p\kern 0.7pt|
\,B^{\kern 1pt 35}$ to $-8\,|\tilde p\kern 0.7pt|\,B^{\kern 1pt 35}$ while $c$ 
runs over the negative interval \mythetag{4.7}. As for the left hand side of
this equation, it is also a continuous function of $c$ that monotonically increases 
from $-160\,|\tilde p\kern 0.7pt|\,B^{\kern 1pt 35}$ to $0$ while $c$ runs over 
this interval. Therefore the equation \mythetag{4.5} has at least one root within 
the negative interval \mythetag{4.7}.\par
    Similarly, if $\tilde p>0$, the right hand side of the equation \mythetag{4.5} 
varies from $8\,|\tilde p\kern 0.7pt|\,B^{\kern 1pt 35}$ to 
$152\,|\tilde p\kern 0.7pt|\,B^{\kern 1pt 35}$ while $c$ runs over the
positive interval \mythetag{4.7} and hence the equation \mythetag{4.5} has at least 
one root within this interval.\par 
     The variable $c$ is related to the original variable $t$ through the formula
\mythetag{4.3}. Therefore from the above considerations we derive the following 
inequalities for $t$:
$$
\gather
\hskip -2em
\gathered
B^{\kern 1pt 2}\,\tilde q^{\kern 0.7pt 6}
+2\,B\,\tilde q^{\kern 0.7pt 4}-2\,B\,\tilde p\,\tilde q^{\kern 1pt 3}
-2\,\tilde q^{\kern 0.7pt 2}-2\,\tilde p\,\tilde q
+\tilde p^{\kern 1pt 2}+\frac{5}{B}\,-\\
-\,\frac{20}{B^{\kern 1pt 2}\,
\tilde q^{\kern 0.7pt 2}}+\frac{10\,\tilde p}{B^{\kern 1pt 2}\,
\tilde q^{\kern 0.7pt 3}}<t<
B^{\kern 1pt 2}\,\tilde q^{\kern 0.7pt 6}
+2\,B\,\tilde q^{\kern 0.7pt 4}-2\,B\,\tilde p
\,\tilde q^{\kern 1pt 3}\,-\\
-\,2\,\tilde q^{\kern 0.7pt 2}-2\,\tilde p\,\tilde q
+\tilde p^{\kern 1pt 2}+\frac{5}{B}
-\frac{20}{B^{\kern 1pt 2}\,
\tilde q^{\kern 0.7pt 2}}
\text{\ \ \ in the case \ \ }\tilde p<0,
\endgathered
\mytag{4.11}\\
\vspace{2ex}
\hskip -2em
\gathered
B^{\kern 1pt 2}\,\tilde q^{\kern 0.7pt 6}
+2\,B\,\tilde q^{\kern 0.7pt 4}-2\,B\,\tilde p\,\tilde q^{\kern 1pt 3}
-2\,\tilde q^{\kern 0.7pt 2}-2\,\tilde p\,\tilde q
+\tilde p^{\kern 1pt 2}+\frac{5}{B}\,-\\
-\,\frac{20}{B^{\kern 1pt 2}\,
\tilde q^{\kern 0.7pt 2}}<t<
B^{\kern 1pt 2}\,\tilde q^{\kern 0.7pt 6}
+2\,B\,\tilde q^{\kern 0.7pt 4}-2\,B\,\tilde p
\,\tilde q^{\kern 1pt 3}-2\,\tilde q^{\kern 0.7pt 2}\,-\\
-\,2\,\tilde p\,\tilde q
+\tilde p^{\kern 1pt 2}+\frac{5}{B}
-\frac{20}{B^{\kern 1pt 2}\,
\tilde q^{\kern 0.7pt 2}}+\frac{10\,\tilde p}{B^{\kern 1pt 2}
\,\tilde q^{\kern 0.7pt 3}}
\text{\ \ \ in the case \ \ }\tilde p>0.
\endgathered
\mytag{4.12}
\endgather
$$
As a result we have proved the following theorem. 
\mytheorem{4.1} If $\tilde p\neq 0$, then for each $\tilde q\geqslant 20
\,\root{3}\of{|\tilde p\kern 0.7pt|}$ and for each value of $B$ in 
\mythetag{4.6} there is at least one real root of the equation \mythetag{2.4} 
satisfying the inequalities \mythetag{4.11} or \mythetag{4.12} respectively.
\endproclaim
     Theorem~\mythetheorem{4.1} means that we have derived the estimate 
\mythetag{4.2} with $C=10\,|\tilde p\kern 0.7pt|/B^{\kern 1pt 2}$ for 
the remainder term $R_1(\tilde p,\tilde q,B)$ in the asymptotic expansion 
\mythetag{4.1} for at least one root of the equation \mythetag{2.4}.\par
     The root $t_2$ in \mythetag{3.20} is handled in a similar way. The 
asymptotic expansion analogous to the expansion \mythetag{4.1} for it is 
written as follows:
$$
\hskip -2em
\aligned
t_2=B^{\kern 1pt 2}\,\tilde q^{\kern 0.7pt 6}-2\,B\,\tilde q^{\kern 0.7pt 4}
&-2\,B\,\tilde p\,\tilde q^{\kern 1pt 3}-2\,\tilde q^{\kern 0.7pt 2}\,+\\
&+\,2\,\tilde p\,\tilde q+\tilde p^{\kern 1pt 2}-\frac{5}{B}
-\frac{20}{B^{\kern 1pt 2}\,\tilde q^{\kern 0.7pt 2}}
+R_2(\tilde p,\tilde q,B).
\endaligned
\mytag{4.13}
$$
The formula \mythetag{4.13} is in agreement with \mythetag{3.19}. 
It means that
$$
\xalignat 4
&\beta_{21}=0,
&&\beta_{22}=-\frac{2}{B},
&&\beta_{23}=-\frac{2\,\tilde p}{B},
&&\beta_{24}=-\frac{2}{B^{\kern 1pt 2}},\\
\displaybreak
&\beta_{25}=\frac{2\,\tilde p}{B^{\kern 1pt 2}},
&&\beta_{26}=\frac{\tilde p^{\kern 1pt 2}}{B^{\kern 1pt 2}}
-\frac{5}{B^{\kern 1pt 3}},
&&\beta_{17}=0
&&\beta_{18}=-\frac{20}{B^{\kern 1pt 4}}
\endxalignat
$$
The formula \mythetag{4.3} in this case is replaced by the following one:
$$
t=B^{\kern 1pt 2}\,\tilde q^{\kern 0.7pt 6}
-2\,B\,\tilde q^{\kern 0.7pt 4}-2\,B\,\tilde p\,\tilde q^{\kern 1pt 3}
-2\,\tilde q^{\kern 0.7pt 2}+2\,\tilde p\,\tilde q
+\tilde p^{\kern 1pt 2}-\frac{5}{B}
-\frac{20}{B^{\kern 1pt 2}\,\tilde q^{\kern 0.7pt 2}}
+\frac{c}{\tilde q^{\kern 0.7pt 3}}.
\quad
\mytag{4.14}
$$
Upon substituting \mythetag{4.14} and \mythetag{4.4} into the equation
\mythetag{2.4} and upon removing denominators the equation \mythetag{2.4} 
is written as an equation very similar to \mythetag{4.5}:
$$
\hskip -2em
16\,c\,B^{\kern 1pt 37}=-80\,\tilde p\,B^{\kern 1pt 35}
+\psi(c,z,\tilde p,B).
\mytag{4.15}
$$
Here $\psi(c,z,\tilde p,B)$ is a polynomial of $c$, $z$, $\tilde p$, and 
$B$ with integer coefficients. The explicit expression for $\psi(c,z,\tilde p,B)$
comprises $1612$ monomials. Therefore it is placed to the ancillary file 
\darkred{{\tt strategy\kern -0.5pt\_\kern 1.5pt formulas\_\kern 0.5pt 04.txt}} 
in a machine-readable form.\par
     Like in the previous case, for each $B$ in \mythetag{4.6} assume that $c$ 
obeys the inequalities
$$
\hskip -2em
\aligned
0<c<\frac{10\,|\tilde p\kern 0.7pt|}{B^{\kern 1pt 2}}\text{\ \ if \ }
\tilde p<0,\\
-\frac{10\,|\tilde p\kern 0.7pt|}{B^{\kern 1pt 2}}<c<0\text{\ \ if \ }
\tilde p>0.
\endaligned
\mytag{4.16}
$$ 
Then assume that $\tilde q$ obeys the inequality \mythetag{4.9}. Under these 
assumptions, taking into account \mythetag{4.8}, one can derive an estimate 
similar to the estimate \mythetag{4.10}:
$$
\hskip -2em
|\psi(c,z,\tilde p,B)|\leqslant 72\,|\tilde p\kern 0.7pt|
\,B^{\kern 1pt 35}.
\mytag{4.17}
$$
Now, writing the inequalities 
$$
\gather
\hskip -2em
\gathered
B^{\kern 1pt 2}\,\tilde q^{\kern 0.7pt 6}
-2\,B\,\tilde q^{\kern 0.7pt 4}-2\,B\,\tilde p\,\tilde q^{\kern 1pt 3}
-2\,\tilde q^{\kern 0.7pt 2}+2\,\tilde p\,\tilde q
+\tilde p^{\kern 1pt 2}-\frac{5}{B}\,-\\
-\,\frac{20}{B^{\kern 1pt 2}\,
\tilde q^{\kern 0.7pt 2}}<t<
B^{\kern 1pt 2}\,\tilde q^{\kern 0.7pt 6}
-2\,B\,\tilde q^{\kern 0.7pt 4}-2\,B\,\tilde p
\,\tilde q^{\kern 1pt 3}-2\,\tilde q^{\kern 0.7pt 2}\,+\\
+\,2\,\tilde p\,\tilde q
+\tilde p^{\kern 1pt 2}-\frac{5}{B}
-\frac{20}{B^{\kern 1pt 2}\,
\tilde q^{\kern 0.7pt 2}}-\frac{10\,\tilde p}{B^{\kern 1pt 2}
\,\tilde q^{\kern 0.7pt 3}}
\text{\ \ \ in the case \ \ }\tilde p<0,
\endgathered
\mytag{4.18}\\
\vspace{2ex}
\hskip -2em
\gathered
B^{\kern 1pt 2}\,\tilde q^{\kern 0.7pt 6}
-2\,B\,\tilde q^{\kern 0.7pt 4}-2\,B\,\tilde p\,\tilde q^{\kern 1pt 3}
-2\,\tilde q^{\kern 0.7pt 2}+2\,\tilde p\,\tilde q
+\tilde p^{\kern 1pt 2}-\frac{5}{B}\,-\\
-\,\frac{20}{B^{\kern 1pt 2}\,
\tilde q^{\kern 0.7pt 2}}-\frac{10\,\tilde p}{B^{\kern 1pt 2}\,
\tilde q^{\kern 0.7pt 3}}<t<
B^{\kern 1pt 2}\,\tilde q^{\kern 0.7pt 6}
-2\,B\,\tilde q^{\kern 0.7pt 4}-2\,B\,\tilde p
\,\tilde q^{\kern 1pt 3}\,-\\
-\,2\,\tilde q^{\kern 0.7pt 2}+2\,\tilde p\,\tilde q
+\tilde p^{\kern 1pt 2}-\frac{5}{B}
-\frac{20}{B^{\kern 1pt 2}\,
\tilde q^{\kern 0.7pt 2}}
\text{\ \ \ in the case \ \ }\tilde p>0,
\endgathered
\mytag{4.19}
\endgather
$$
and then applying the estimate \mythetag{4.17} to the equation \mythetag{4.15} 
and taking into account \mythetag{4.16}, we can easily prove the following
theorem. 
\mytheorem{4.2} If $\tilde p\neq 0$, then for each $\tilde q\geqslant 20
\,\root{3}\of{|\tilde p\kern 0.7pt|}$ and for each value of $B$ in 
\mythetag{4.6} there is at least one real root of the equation \mythetag{2.4} 
satisfying the inequalities \mythetag{4.18} or \mythetag{4.19} respectively.
\endproclaim
     Theorem~\mythetheorem{4.2} is similar to Theorem~\mythetheorem{4.1}.
It means that we have got the estimate 
$$
\hskip -2em
|R_2(\tilde p,\tilde q,B)|<\frac{C}{\tilde q^{\kern 0.7pt 3}}
\text{, \ where \ }C=\frac{10\,|\tilde p\kern 0.7pt|}{B^{\kern 1pt 2}},
\mytag{4.20}
$$
for the remainder term $R_2(\tilde p,\tilde q,B)$ of the
asymptotic expansion \mythetag{4.13}.\par
     The root $t_3$ in \mythetag{3.20} is somewhat different from 
$t_1$ and $t_2$. The analog of the asymp\-totic expansions \mythetag{4.1}
and \mythetag{4.13} for this root is written as 
$$
\hskip -2em
t_3=B\,\tilde q^{\kern 0.7pt 4}-\tilde p\,\tilde q+\frac{16}{B}
+R_3(\tilde p,\tilde q,B).
\mytag{4.21}
$$
Like in \mythetag{4.2} and \mythetag{4.20}, in this case we shall derive
the inverse cubic estimate 
$$
\hskip -2em
|R_3(\tilde p,\tilde q,B)|<\frac{C}{\tilde q^{\kern 0.7pt 3}}
\mytag{4.22}
$$
for the remainder term in \mythetag{4.21}. For this purpose we substitute 
$$
\hskip -2em
t=B\,\tilde q^{\kern 0.7pt 4}-\tilde p\,\tilde q+\frac{16}{B}
+\frac{c}{\tilde q^{\kern 0.7pt 3}}
\mytag{4.23}
$$
into the equation \mythetag{2.4}. Immediately after that we perform the
substitution \mythetag{4.4} into the equation obtained by substituting 
\mythetag{4.23} into \mythetag{2.4}. Upon removing denominators, the 
resulting equation can be written in the following form:
$$
\hskip -2em
2\,c\,B^{\kern 1pt 23}=32\,\tilde p\,B^{\kern 1pt 21}+f(c,z,\tilde p,B).
\mytag{4.24}
$$
Here $f(c,z,\tilde p,B)$ is a polynomial of $c$, $z$, $\tilde p$, and 
$B$ with integer coefficients. The explicit expression for 
$f(c,z,\tilde p,B)$ comprises $490$ monomials. Therefore it is placed to 
the ancillary file \darkred{{\tt strategy\kern -0.5pt\_\kern 1.5pt 
formulas\_\kern 0.5pt 04.txt}} in a machine-readable form.\par
For each $B$ in \mythetag{4.6} assume that $c$ obeys the condition
$$
\hskip -2em
\aligned
-\frac{32\,|\tilde p\kern 0.7pt|}{B^{\kern 1pt 2}}<c<0\text{\ \ if \ }
\tilde p<0,\\
0<c<\frac{32\,|\tilde p\kern 0.7pt|}{B^{\kern 1pt 2}}\text{\ \ if \ }
\tilde p>0.
\endaligned
\mytag{4.25}
$$ 
Again, the case $\tilde p=0$ is exceptional. It should be studied separately.
\par
     Like in \mythetag{4.9}, assume that $\tilde q$ obeys the inequality
$$
\hskip -2em
\tilde q\geqslant 7\,\root{3}\of{|\tilde p\kern 0.7pt|}.
\mytag{4.26}
$$
Let's apply \mythetag{4.26} to \mythetag{4.4}. As  a result we derive the 
inequality $|z|\leqslant 1/7\,|\tilde p\kern 0.7pt|^{-1/3}$. 
Applying this inequality along with the inequalities \mythetag{4.25} and 
\mythetag{4.8} to the polynomial $f(c,z,\tilde p,B)$ in \mythetag{4.24}, 
we derive the following estimate for it: 
$$
\hskip -2em
|f(c,z,\tilde p,B)|\leqslant 26\,|\tilde p\kern 0.7pt|
\,B^{\kern 1pt 21}.
\mytag{4.27}
$$
The inequality \mythetag{4.27} holds for each value of $B$ in \mythetag{4.6}.
\par
     Now we apply the inequality \mythetag{4.27} to the equation \mythetag{4.24}.
If $\tilde p<0$, it means that the right hand side of the equation \mythetag{4.24} 
is a continuous function of $c$ that varies from $-58\,|\tilde p\kern 0.7pt|
\,B^{\kern 1pt 35}$ to $-6\,|\tilde p\kern 0.7pt|\,B^{\kern 1pt 35}$ while $c$ 
runs over the negative interval \mythetag{4.25}. As for the left hand side of
this equation, it is also a continuous function of $c$ that monotonically increases 
from $-64\,|\tilde p\kern 0.7pt|\,B^{\kern 1pt 35}$ to $0$ while $c$ runs over 
this interval. Therefore the equation \mythetag{4.24} has at least one root within 
the negative interval \mythetag{4.25}. Similarly, if $\tilde p>0$, the right hand 
side of the equation \mythetag{4.24} varies from $6\,|\tilde p\kern 0.7pt|
\,B^{\kern 1pt 35}$ to $58\,|\tilde p\kern 0.7pt|\,B^{\kern 1pt 35}$ while $c$ runs 
over the positive interval \mythetag{4.25} and hence the equation \mythetag{4.24} 
has at least one root within this interval.\par 
     The variable $c$ is related to the original variable $t$ through the formula
\mythetag{4.23}. Therefore from the above considerations we derive the following 
inequalities for $t$:
$$
\gather
\hskip -2em
B\,\tilde q^{\kern 0.7pt 4}-\tilde p\,\tilde q+\frac{16}{B}
+\frac{32\,\tilde p}{B^{\kern 1pt 2}\,
\tilde q^{\kern 0.7pt 3}}<t<
B\,\tilde q^{\kern 0.7pt 4}-\tilde p\,\tilde q+\frac{16}{B}
\text{\ \ \ if \ \ }\tilde p<0,
\mytag{4.28}\\
\vspace{1ex}
\hskip -2em
B\,\tilde q^{\kern 0.7pt 4}-\tilde p\,\tilde q+\frac{16}{B}
<t<B\,\tilde q^{\kern 0.7pt 4}-\tilde p\,\tilde q+\frac{16}{B}
+\frac{32\,\tilde p}{B^{\kern 1pt 2}\,
\tilde q^{\kern 0.7pt 3}}
\text{\ \ \ if \ \ }\tilde p>0.
\mytag{4.29}
\endgather
$$
As a result we have proved the following theorem. 
\mytheorem{4.3} If $\tilde p\neq 0$, then for each $\tilde q\geqslant 7
\,\root{3}\of{|\tilde p\kern 0.7pt|}$ and for each value of $B$ in 
\mythetag{4.6} there is at least one real root of the equation \mythetag{2.4} 
satisfying the inequalities \mythetag{4.28} or \mythetag{4.29} respectively.
\endproclaim
     Theorem~\mythetheorem{4.3} means that we have derived the estimate 
\mythetag{4.22} with $C=32\,|\tilde p\kern 0.7pt|/B^{\kern 1pt 2}$ for 
the remainder term $R_3(\tilde p,\tilde q,B)$ in the asymptotic expansion 
\mythetag{4.21} for at least one root of the equation \mythetag{2.4}.\par
\head
5. Asymptotic estimates for complex roots.
\endhead
     There are two complex roots $t_4$ and $t_5$ of the equation 
\mythetag{2.4} satisfying the condition \mythetag{2.5}. The leading terms
of their asymptotics as $\tilde q\to +\infty$ are given by the formulas
\mythetag{3.21}. Specifying the first formula \mythetag{3.21}, we write
$$
\hskip -2em
t_4=(\sqrt{2}+1)\,\goth i\,\tilde q^{\kern 0.7pt 2}
-\frac{(10+7\,\sqrt{2}\,)\,\goth i\,}{B^{\kern 1pt 2}\,\tilde q^{\kern 0.7pt 2}}
+R_4(\tilde p,\tilde q,B).
\mytag{5.1}
$$
Here $\goth i=\sqrt{-1}$. Our goal is to derive an estimate of the form 
$$
\hskip -2em
|R_4(\tilde p,\tilde q,B)|<\frac{C}{\tilde q^{\kern 0.7pt 5}}
\mytag{5.2}
$$
for the remainder term in \mythetag{5.1}. In order to get such an estimate 
we substitute 
$$
\hskip -2em
t=(\sqrt{2}+1)\,\goth i\,\tilde q^{\kern 0.7pt 2}
-\frac{(10+7\,\sqrt{2}\,)\,\goth i\,}{B^{\kern 1pt 2}
\,\tilde q^{\kern 0.7pt 2}}
+\frac{c\,\goth i}{\tilde q^{\kern 0.7pt 5}}
\mytag{5.3}
$$
into the equation \mythetag{2.4}. Immediately after that we perform the
substitution \mythetag{4.4} into the equation obtained by substituting 
\mythetag{5.3} into \mythetag{2.4}. Upon removing denominators, the 
resulting equation can be written in the following form:
$$
\hskip -2em
16\,c\,B^{\kern 1pt 30}=-32\,(10+7\,\sqrt{2})\,\tilde p\,B^{\kern 1pt 27}
+\eta(c,z,\tilde p,B).
\mytag{5.4}
$$
Here $\eta(c,z,\tilde p,B)$ is a polynomial of $c$, $z$, $\tilde p$, and 
$B$. The fully expanded expression for $\eta(c,z,\tilde p,B)$ comprises 
$1031$ monomials. Therefore it is placed to the ancillary file 
\darkred{{\tt strategy\kern -0.5pt\_\kern 1.5pt formulas\_\kern 0.5pt 04.txt}} 
in a machine-readable form.\par
     Note that the value of the irrational coefficient in \mythetag{5.4} obeys 
the inequalities
$$
\hskip -2em
636<32\,(10+7\,\sqrt{2})\approx 636.78<640=16\cdot 40.
\mytag{5.5}
$$
Therefore for each $B$ in \mythetag{4.6} assume that $c$ obeys the condition
$$
\hskip -2em
\aligned
0<c<\frac{80\,|\tilde p\kern 0.7pt|}{B^{\kern 1pt 3}}\text{\ \ if \ }
\tilde p<0,\\
-\frac{80\,|\tilde p\kern 0.7pt|}{B^{\kern 1pt 3}}<c<0\text{\ \ if \ }
\tilde p>0.
\endaligned
\mytag{5.6}
$$ 
Like in \mythetag{4.7}, \mythetag{4.16}, and \mythetag{4.25} above, the 
case $\tilde p=0$ is exceptional. It should be studied separately.\par 
     Assume that $\tilde q$ obeys the following inequality similar to \mythetag{4.9} 
and \mythetag{4.26}:
$$
\hskip -2em
\tilde q\geqslant 15\,\root{3}\of{|\tilde p\kern 0.7pt|}.
\mytag{5.7}
$$
Let's apply \mythetag{5.7} to \mythetag{4.4}. As  a result we derive the 
inequality $|z|\leqslant 1/15\,|\tilde p\kern 0.7pt|^{-1/3}$. Applying this 
inequality along with the inequalities \mythetag{5.5} and 
\mythetag{4.8} to the polynomial $\eta(c,z,\tilde p,B)$ in \mythetag{5.4}, 
we derive the following estimate for it: 
$$
\hskip -2em
|\eta(c,z,\tilde p,B)|\leqslant 512\,|\tilde p\kern 0.7pt|
\,B^{\kern 1pt 27}.
\mytag{5.8}
$$\par 
     Let's compare \mythetag{5.8} with \mythetag{5.5} and then apply 
the inequality \mythetag{5.8} to the equation \mythetag{5.4}. If $\tilde p<0$, 
it means that the right hand side of the equation \mythetag{5.4} 
is a continuous function of $c$ that varies from $124\,|\tilde p\kern 0.7pt|
\,B^{\kern 1pt 27}$ to $1152\,|\tilde p\kern 0.7pt|\,B^{\kern 1pt 27}$ 
while $c$ runs over the positive interval \mythetag{5.6}. As for the left hand 
side of this equation, it is also a continuous function of $c$ that monotonically 
increases from $0$ to $1280\,|\tilde p\kern 0.7pt|\,B^{\kern 1pt 27}$ while $c$ 
runs over this interval. Therefore the equation \mythetag{5.4} has at least one 
root within the positive interval \mythetag{5.6}. Similarly, if $\tilde p>0$, 
the right hand side of the equation \mythetag{5.4} varies from $-1152\,|\tilde p
\kern 0.7pt|\,B^{\kern 1pt 27}$ to $-124\,|\tilde p\kern 0.7pt|\,B^{\kern 1pt 27}$ 
while $c$ runs over the negative interval \mythetag{5.6} and hence the equation 
\mythetag{5.4} has at least one root within this interval.\par 
     The variable $c$ is related to the original variable $t$ through the formula
\mythetag{5.3}. Therefore from the above considerations we derive the following 
inequalities for $t$:
$$
\gather
\hskip -2em
\gathered
(\sqrt{2}+1)\,\tilde q^{\kern 0.7pt 2}
-\frac{10+7\,\sqrt{2}\,}{B^{\kern 1pt 2}\,\tilde q^{\kern 0.7pt 2}}
<\Img t<(\sqrt{2}+1)\,\tilde q^{\kern 0.7pt 2}\,-\\
-\,\frac{10+7\,\sqrt{2}\,}{B^{\kern 1pt 2}\,\tilde q^{\kern 0.7pt 2}}
-\frac{80\,\tilde p}{B^{\kern 1pt 3}\,\tilde q^{\kern 0.7pt 5}}
\text{\ \ \ in the case \ \ }\tilde p<0,
\endgathered
\mytag{5.9}\\
\vspace{1ex}
\hskip -2em
\gathered
(\sqrt{2}+1)\,\tilde q^{\kern 0.7pt 2}
-\frac{10+7\,\sqrt{2}\,}{B^{\kern 1pt 2}\,\tilde q^{\kern 0.7pt 2}}
-\frac{80\,\tilde p}{B^{\kern 1pt 3}\,\tilde q^{\kern 0.7pt 5}}
<\Img t<(\sqrt{2}+1)\,\tilde q^{\kern 0.7pt 2}\,-\\
-\,\frac{10+7\,\sqrt{2}\,}{B^{\kern 1pt 2}\,\tilde q^{\kern 0.7pt 2}}
\text{\ \ \ in the case \ \ }\tilde p>0. 
\endgathered
\mytag{5.10}
\endgather
$$
The inequalities \mythetag{5.9} and \mythetag{5.10} are analogs of 
the corresponding inequalities in the case of the real roots, \pagebreak 
see \mythetag{4.11} and \mythetag{4.12}, \mythetag{4.18} and \mythetag{4.19}, 
or \mythetag{4.28} and \mythetag{4.29}. They mean that we have proved the 
following theorem.    
\mytheorem{5.1} If $\tilde p\neq 0$, then for each $\tilde q\geqslant 15
\,\root{3}\of{|\tilde p\kern 0.7pt|}$ and for each value of $B$ in 
\mythetag{4.6} there is at least one purely imaginary root of the equation 
\mythetag{2.4} satisfying the inequalities \mythetag{5.9} or \mythetag{5.10} 
respectively.
\endproclaim
     Theorem~\mythetheorem{5.1} means that we have derived the estimate 
\mythetag{5.2} with $C=80\,|\tilde p\kern 0.7pt|/B^{\kern 1pt 3}$ for 
the remainder term $R_4(\tilde p,\tilde q,B)$ in the asymptotic expansion 
\mythetag{5.1} for at least one root of the equation \mythetag{2.4}.\par
     The complex root $t_5$ is similar to the root $t_4$. For this root we 
write
$$
\hskip -2em
t_5=(\sqrt{2}-1)\,\goth i\,\tilde q^{\kern 0.7pt 2}
+\frac{(10-7\,\sqrt{2}\,)\,\goth i\,}{B^{\kern 1pt 2}\,\tilde q^{\kern 0.7pt 2}}
+R_5(\tilde p,\tilde q,B).
\mytag{5.11}
$$
Here $\goth i=\sqrt{-1}$. Like in \mythetag{5.1}, our goal is to derive an 
estimate of the form 
$$
\hskip -2em
|R_5(\tilde p,\tilde q,B)|<\frac{C}{\tilde q^{\kern 0.7pt 5}}
\mytag{5.12}
$$
for the remainder term in \mythetag{5.11}. In order to get such an estimate 
we substitute 
$$
\hskip -2em
t=(\sqrt{2}-1)\,\goth i\,\tilde q^{\kern 0.7pt 2}
+\frac{(10-7\,\sqrt{2}\,)\,\goth i\,}{B^{\kern 1pt 2}
\,\tilde q^{\kern 0.7pt 2}}
+\frac{c\,\goth i}{\tilde q^{\kern 0.7pt 5}}
\mytag{5.13}
$$
into the equation \mythetag{2.4}. Immediately after that we perform the
substitution \mythetag{4.4} into the equation obtained by substituting 
\mythetag{5.13} into \mythetag{2.4}. Upon removing denominators, the 
resulting equation can be written in the following form:
$$
\hskip -2em
16\,c\,B^{\kern 1pt 30}=32\,(10-7\,\sqrt{2})\,\tilde p\,B^{\kern 1pt 27}
+\zeta(c,z,\tilde p,B).
\mytag{5.14}
$$
Here $\zeta(c,z,\tilde p,B)$ is a polynomial of $c$, $z$, $\tilde p$, and 
$B$. The fully expanded expression for $\zeta(c,z,\tilde p,B)$ comprises 
$1031$ monomials. Therefore it is placed to the ancillary file 
\darkred{{\tt strategy\kern -0.5pt\_\kern 1.5pt formulas\_\kern 0.5pt 04.txt}} 
in a machine-readable form.\par
     Note that the value of the irrational coefficient in \mythetag{5.14} obeys 
the inequality
$$
\hskip -2em
3<32\,(10-7\,\sqrt{2})\approx 3.22<4=16\kern 0.8pt\cdot\frac{1}{4}.
\mytag{5.15}
$$
The number in \mythetag{5.15} is substantially smaller than the number in
\mythetag{5.5}. Therefore, instead of \mythetag{5.6} we write the following 
inequalities:
$$
\hskip -2em
\aligned
-\frac{|\tilde p\kern 0.7pt|}{2\,B^{\kern 1pt 3}}<c<0\text{\ \ if \ }
\tilde p<0,\\
0<c<\frac{|\tilde p\kern 0.7pt|}{2\,B^{\kern 1pt 3}}\text{\ \ if \ }
\tilde p>0.
\endaligned
\mytag{5.16}
$$ 
The inequalities \mythetag{5.16} should be fulfilled for each value of $B$ 
in \mythetag{4.6}.\par
     Assume that $\tilde q$ obeys the following inequality similar to 
\mythetag{4.9}, \mythetag{4.26}, and \mythetag{5.7}:
$$
\hskip -2em
\tilde q\geqslant 3600\,\root{3}\of{|\tilde p\kern 0.7pt|}.
\mytag{5.17}
$$
The coefficient $3600$ in \mythetag{5.17} is much greater than the 
coefficient $15$ in \mythetag{5.7}. Applying \mythetag{5.17} to 
\mythetag{4.4}, we derive  $|z|\leqslant 1/3600\,|\tilde p\kern 0.7pt|^{-1/3}$.
Applying this inequality along with the inequalities \mythetag{5.15} and 
\mythetag{4.8} to the polynomial $\zeta(c,z,\tilde p,B)$ in the equation
\mythetag{5.14}, we derive the following estimate for it: 
$$
\hskip -2em
|\zeta(c,z,\tilde p,B)|\leqslant 2\,|\tilde p\kern 0.7pt|
\,B^{\kern 1pt 27}.
\mytag{5.18}
$$\par 
     Using the estimate \mythetag{5.18}, we can easily prove that 
under the assumption \mythetag{5.17} the equation \mythetag{5.14} has at 
least one root obeying the corresponding inequalities \mythetag{5.16}. 
The variable $c$ is related to the original variable $t$ through the formula
\mythetag{5.13}. Therefore from \mythetag{5.16} we derive the following 
inequalities for $t$:
$$
\gather
\hskip -2em
\gathered
(\sqrt{2}-1)\,\tilde q^{\kern 0.7pt 2}
+\frac{10-7\,\sqrt{2}\,}{B^{\kern 1pt 2}\,\tilde q^{\kern 0.7pt 2}}
+\frac{\tilde p}{2\,B^{\kern 1pt 3}\,\tilde q^{\kern 0.7pt 5}}
<\Img t<(\sqrt{2}-1)\,\tilde q^{\kern 0.7pt 2}\,+\\
+\,\frac{10-7\,\sqrt{2}\,}{B^{\kern 1pt 2}\,\tilde q^{\kern 0.7pt 2}}
\text{\ \ \ in the case \ \ }\tilde p<0,
\endgathered
\mytag{5.19}\\
\vspace{1ex}
\hskip -2em
\gathered
(\sqrt{2}-1)\,\tilde q^{\kern 0.7pt 2}
+\frac{10-7\,\sqrt{2}\,}{B^{\kern 1pt 2}\,\tilde q^{\kern 0.7pt 2}}
<\Img t<(\sqrt{2}-1)\,\tilde q^{\kern 0.7pt 2}\,+\\
+\,\frac{10-7\,\sqrt{2}\,}{B^{\kern 1pt 2}\,\tilde q^{\kern 0.7pt 2}}
+\frac{\tilde p}{2\,B^{\kern 1pt 3}\,\tilde q^{\kern 0.7pt 5}}
\text{\ \ \ in the case \ \ }\tilde p>0. 
\endgathered
\mytag{5.20}
\endgather
$$
As a result we have proved the following theorem. 
\mytheorem{5.2} If $\tilde p\neq 0$, then for each $\tilde q\geqslant 3600
\,\root{3}\of{|\tilde p\kern 0.7pt|}$ and for each value of $B$ in 
\mythetag{4.6} there is at least one purely imaginary root of the equation 
\mythetag{2.4} satisfying the inequalities \mythetag{5.19} or \mythetag{5.20} 
respectively.
\endproclaim
     Theorem~\mythetheorem{5.2} means that we have derived the estimate 
\mythetag{5.12} with the coefficient $C=|\tilde p\kern 0.7pt|/(2
\,B^{\kern 1pt 3})$ for the remainder term $R_5(\tilde p,\tilde q,B)$ in 
the asymptotic expansion \mythetag{5.11} for at least one root of the 
equation \mythetag{2.4}.\par
\head
6. Non-intersection of asymptotic sites.
\endhead
    In the previous two sections we have used four inequalities \mythetag{4.9},
\mythetag{4.26}, \mythetag{5.7}, and \mythetag{5.17}. The inequality 
\mythetag{5.17} is the strongest of them. It looks like 
$$
\hskip -2em
\tilde q\geqslant 3600\,\root{3}\of{|\tilde p\kern 0.7pt|}.
\mytag{6.1}
$$
If the inequality \mythetag{6.1} is fulfilled, then all of the inequalities 
\mythetag{4.9}, \mythetag{4.26}, \mythetag{5.7}, \mythetag{5.17} are also
fulfilled. Due to these inequalities, applying Theorems~\mythetheorem{4.1},
\mythetheorem{4.2}, \mythetheorem{4.3}, \mythetheorem{5.1}, and 
\mythetheorem{5.2}, we find that there are five asymptotic sites, each of
which comprises at least one root of the equation \mythetag{2.4}.\par 
     The site comprising the root $t_3$ is given by the inequalities
\mythetag{4.28} and \mythetag{4.29}. From \mythetag{6.1} and \mythetag{4.8}
we derive the following inequalities:
$$
\xalignat 2
&\hskip -2em
\tilde q^{\kern 0.7pt 3}\geqslant 3600^{\kern 0.7pt 3}\,|\tilde p\kern 0.7pt|
&&\tilde q\geqslant 3600,\\
\vspace{-1.5ex}
\mytag{6.2}\\
\vspace{-1.5ex}
&\hskip -2em
\frac{\tilde q^{\kern 0.7pt 3}}{|\tilde p\kern 0.7pt|}
\geqslant 3600^{\kern 0.7pt 3},
&&\frac{|\tilde p\kern 0.7pt|}{\tilde q^{\kern 0.7pt 3}}
\leqslant\frac{1}{3600^{\kern 0.7pt 3}}.
\endxalignat
$$
From \mythetag{4.6} we easily derive the following inequalities for $B$:
$$
\xalignat 2
&\hskip -2em
1\leqslant B\leqslant 9,
&&\frac{1}{9}\leqslant \frac{1}{B}\leqslant 1.
\mytag{6.3}
\endxalignat
$$
Applying \mythetag{6.2} and \mythetag{6.3} to the left hand side 
of \mythetag{4.28}, we derive
$$
\hskip -2em
\gathered
B\,\tilde q^{\kern 0.7pt 4}-\tilde p\,\tilde q+\frac{16}{B}
+\frac{32\,\tilde p}{B^{\kern 1pt 2}\,
\tilde q^{\kern 0.7pt 3}}\geqslant B\,\tilde q^{\kern 0.7pt 4}
\Bigl(1-\frac{1}{3600^{\kern 0.7pt 3}}\Bigr)+\frac{16}{B}
-\frac{32}{B^{\kern 1pt 2}\,
3600^{\kern 0.7pt 3}}\,\geqslant\\
\geqslant\,3600^{\kern 0.7pt 4}\Bigl(1
-\frac{1}{3600^{\kern 0.7pt 3}}\Bigr)+\frac{16}{9}
-\frac{32}{3600^{\kern 0.7pt 3}}
\approx 1.68\cdot10^{\kern 0.5pt 14}>0.
\endgathered
\mytag{6.4}
$$
The left hand side of \mythetag{4.29} is treated similarly:
$$
\hskip -2em
\gathered
B\,\tilde q^{\kern 0.7pt 4}-\tilde p\,\tilde q+\frac{16}{B}
\geqslant B\,\tilde q^{\kern 0.7pt 4}
\Bigl(1-\frac{1}{3600^{\kern 0.7pt 3}}\Bigr)+\frac{16}{B}
\,\geqslant\\
\geqslant\,3600^{\kern 0.7pt 4}\Bigl(1
-\frac{1}{3600^{\kern 0.7pt 3}}\Bigr)+\frac{16}{9}
\approx 1.67\cdot 10^{\kern 0.5pt 14}>0.
\endgathered
\mytag{6.5}
$$\par
     The root $t_2$ is delimited by the inequalities \mythetag{4.18}
and \mythetag{4.19}. Let's compare the left hand side of \mythetag{4.18}
with the right hand side of \mythetag{4.28}. For their difference, using
\mythetag{6.2} and \mythetag{6.3}, we derive the following estimate:
$$
\gathered
\Bigl(B^{\kern 1pt 2}\,\tilde q^{\kern 0.7pt 6}
-2\,B\,\tilde q^{\kern 0.7pt 4}-2\,B\,\tilde p\,\tilde q^{\kern 1pt 3}
-2\,\tilde q^{\kern 0.7pt 2}+2\,\tilde p\,\tilde q
+\tilde p^{\kern 1pt 2}-\frac{5}{B}-\frac{20}{B^{\kern 1pt 2}\,
\tilde q^{\kern 0.7pt 2}}\Bigr)-\\
-\Bigl(B\,\tilde q^{\kern 0.7pt 4}-\tilde p\,\tilde q+\frac{16}{B}
\Bigr)=B^{\kern 1pt 2}\,\tilde q^{\kern 0.7pt 6}\Bigl(1-\frac{3}
{B\,\tilde q^{\kern 0.7pt 2}}
-\frac{2\,\tilde p}{B\,\tilde q^{\kern 1pt 3}}
-\frac{2}{B^{\kern 1pt 2}\,\tilde q^{\kern 0.7pt 4}}
+\frac{3\,\tilde p}{B^{\kern 1pt 2}\,\tilde q^{\kern 0.7pt 5}}\,-\\
-\,\frac{21}{B^{\kern 1pt 3}\,\tilde q^{\kern 0.7pt 6}}
-\frac{20}{B^{\kern 1pt 4}\,\tilde q^{\kern 0.7pt 8}}\Bigr)
\geqslant B^{\kern 1pt 2}\,\tilde q^{\kern 0.7pt 6}\Bigl(1
-\frac{3}{3600^{\kern 0.7pt 2}}-\frac{2}{3600^{\kern 0.7pt 3}}
-\frac{2}{3600^{\kern 0.7pt 4}}-\frac{3}{3600^{\kern 0.7pt 5}}\,-\\
-\frac{21}{3600^{\kern 0.7pt 6}}-\frac{20}{3600^{\kern 0.7pt 8}}
\Bigr)\approx 0.99\cdot B^{\kern 1pt 2}\,\tilde q^{\kern 0.7pt 6}
\geqslant 0.99\cdot 3600^{\kern 0.7pt 6}\approx 2.18\cdot
10^{\kern 0.5pt 21}>0.
\endgathered
\quad
\mytag{6.6}
$$
The difference of the left hand side of \mythetag{4.19} and
the right hand side of \mythetag{4.29} is treated similarly.
For this difference we have 
$$
\gathered
\Bigl(B^{\kern 1pt 2}\,\tilde q^{\kern 0.7pt 6}
-2\,B\,\tilde q^{\kern 0.7pt 4}-2\,B\,\tilde p\,\tilde q^{\kern 1pt 3}
-2\,\tilde q^{\kern 0.7pt 2}+2\,\tilde p\,\tilde q
+\tilde p^{\kern 1pt 2}-\frac{5}{B}-\frac{20}{B^{\kern 1pt 2}\,
\tilde q^{\kern 0.7pt 2}}\,-\\
-\,\frac{10\,\tilde p}{B^{\kern 1pt 2}\,\tilde q^{\kern 0.7pt 3}}
\Bigr)-\Bigl(B\,\tilde q^{\kern 0.7pt 4}-\tilde p\,\tilde q+\frac{16}{B}
+\frac{32\,\tilde p}{B^{\kern 1pt 2}\,
\tilde q^{\kern 0.7pt 3}}\Bigr)\geqslant B^{\kern 1pt 2}
\,\tilde q^{\kern 0.7pt 6}\Bigl(1-\frac{3}{3600^{\kern 0.7pt 2}}\,-\\
-\,\frac{2}{3600^{\kern 0.7pt 3}}-\frac{2}{3600^{\kern 0.7pt 4}}
-\frac{3}{3600^{\kern 0.7pt 5}}-\frac{21}{3600^{\kern 0.7pt 6}}
-\frac{20}{3600^{\kern 0.7pt 8}}-\frac{42}{3600^{\kern 0.7pt 9}}
\Bigr)\approx\\
\vspace{2ex}
\approx 0.99\cdot B^{\kern 1pt 2}\,\tilde q^{\kern 0.7pt 6}
\geqslant 0.99\cdot 3600^{\kern 0.7pt 6}\approx 2.17\cdot
10^{\kern 0.5pt 21}>0.
\endgathered
\quad
\mytag{6.7}
$$\par 
     The root $t_1$ is delimited by the inequalities \mythetag{4.11}
and \mythetag{4.12}. Let's consider the difference of the left hand side of 
\mythetag{4.11} and the right hand side of \mythetag{4.18}:
$$
\gather
\Bigl(B^{\kern 1pt 2}\,\tilde q^{\kern 0.7pt 6}
+2\,B\,\tilde q^{\kern 0.7pt 4}-2\,B\,\tilde p\,\tilde q^{\kern 1pt 3}
-2\,\tilde q^{\kern 0.7pt 2}-2\,\tilde p\,\tilde q
+\tilde p^{\kern 1pt 2}+\frac{5}{B}-\frac{20}{B^{\kern 1pt 2}\,
\tilde q^{\kern 0.7pt 2}}\,+\\
+\,\frac{10\,\tilde p}{B^{\kern 1pt 2}\,
\tilde q^{\kern 0.7pt 3}}
\Bigr)-\Bigl(B^{\kern 1pt 2}\,\tilde q^{\kern 0.7pt 6}
-2\,B\,\tilde q^{\kern 0.7pt 4}-2\,B\,\tilde p
\,\tilde q^{\kern 1pt 3}-2\,\tilde q^{\kern 0.7pt 2}
+2\,\tilde p\,\tilde q+\tilde p^{\kern 1pt 2}\,-\\
\displaybreak
-\,\frac{5}{B}-\frac{20}{B^{\kern 1pt 2}\,\tilde q^{\kern 0.7pt 2}}
-\frac{10\,\tilde p}{B^{\kern 1pt 2}\,\tilde q^{\kern 0.7pt 3}}\Bigr)
=4\,B\,\tilde q^{\kern 0.7pt 4}-4\,\tilde p\,\tilde q+\frac{10}{B}
+\frac{20\,\tilde p}{B^{\kern 1pt 2}\,\tilde q^{\kern 0.7pt 3}}.
\endgather
$$
Applying \mythetag{6.2} and \mythetag{6.3} to the above expression, 
we get the following estimate:
$$
\hskip -2em
\gathered
4\,B\,\tilde q^{\kern 0.7pt 4}-4\,\tilde p\,\tilde q+\frac{10}{B}
+\frac{20\,\tilde p}{B^{\kern 1pt 2}\,\tilde q^{\kern 0.7pt 3}}
=4\,B\,\tilde q^{\kern 0.7pt 4}\Bigl(1-\frac{\tilde p}{B\,
\tilde q^{\kern 0.7pt 3}}+\frac{5}{2\,B^{\kern 1pt 2}
\,\tilde q^{\kern 0.7pt 4}}\,+\\
+\,\frac{5\,\tilde p}{B^{\kern 1pt 3}\,\tilde q^{\kern 0.7pt 7}}
\Bigr)\geqslant 4\,B\,\tilde q^{\kern 0.7pt 4}\Bigl(1
-\frac{1}{3600^{\kern 0.7pt 3}}-\frac{5}{2\cdot 3600^{\kern 0.7pt 4}}
-\frac{5}{3600^{\kern 0.7pt 7}}\Bigr)\,\approx\\
\vspace{2ex}
\approx\,4.00\cdot B\,\tilde q^{\kern 0.7pt 4}
\geqslant 4.00\cdot 3600^{\kern 0.7pt 4}\approx 6.71\cdot
10^{\kern 0.5pt 14}>0.
\endgathered
\mytag{6.8}
$$ 
The difference of the left hand side of \mythetag{4.12} and
the right hand side of \mythetag{4.19} is treated similarly.
For this difference we have 
$$
\gathered
\Bigl(B^{\kern 1pt 2}\,\tilde q^{\kern 0.7pt 6}
+2\,B\,\tilde q^{\kern 0.7pt 4}-2\,B\,\tilde p\,\tilde q^{\kern 1pt 3}
-2\,\tilde q^{\kern 0.7pt 2}-2\,\tilde p\,\tilde q
+\tilde p^{\kern 1pt 2}+\frac{5}{B}-\frac{20}{B^{\kern 1pt 2}
\,\tilde q^{\kern 0.7pt 2}}\Bigr)\,-\\
-\Bigl(B^{\kern 1pt 2}\,\tilde q^{\kern 0.7pt 6}
-2\,B\,\tilde q^{\kern 0.7pt 4}-2\,B\,\tilde p
\,\tilde q^{\kern 1pt 3}-2\,\tilde q^{\kern 0.7pt 2}
+2\,\tilde p\,\tilde q+\tilde p^{\kern 1pt 2}-\frac{5}{B}
-\frac{20}{B^{\kern 1pt 2}\,\tilde q^{\kern 0.7pt 2}}\Bigr)=\\
=4\,B\,\tilde q^{\kern 0.7pt 4}-4\,\tilde p\,\tilde q+\frac{10}{B}
\geqslant 4\,B\,\tilde q^{\kern 0.7pt 4}
\Bigl(1-\frac{1}{3600^{\kern 0.7pt 3}}
-\frac{5}{2\cdot 3600^{\kern 0.7pt 4}}\Bigr)\,\approx\\
\vspace{2ex}
\approx\,4.00\cdot B\,\tilde q^{\kern 0.7pt 4}
\geqslant 4.00\cdot 3600^{\kern 0.7pt 4}\approx 6.71\cdot
10^{\kern 0.5pt 14}>0.
\endgathered
\mytag{6.9}
$$
The inequalities \mythetag{6.4}, \mythetag{6.5}, \mythetag{6.6}, 
\mythetag{6.7}, \mythetag{6.8}, and \mythetag{6.9} mean that if the
inequality \mythetag{6.1} is fulfilled, then the asymptotic sites 
for the roots $t_1$, $t_2$, $t_3$ do not intersect with each other
and are located within the positive half-line of the real axis.
\par 
     Now let's proceed to the purely imaginary roots $t_4$ and $t_5$.
The root $t_5$ is delimited by the inequalities \mythetag{5.19} and
\mythetag{5.20}. Applying \mythetag{6.2} and \mythetag{6.3} to the 
left hand side of \mythetag{5.19}, we derive the following estimate
for it:
$$
\hskip -2em
\gathered
(\sqrt{2}-1)\,\tilde q^{\kern 0.7pt 2}
+\frac{10-7\,\sqrt{2}\,}{B^{\kern 1pt 2}\,\tilde q^{\kern 0.7pt 2}}
+\frac{\tilde p}{2\,B^{\kern 1pt 3}\,\tilde q^{\kern 0.7pt 5}}
=(\sqrt{2}-1)\,\tilde q^{\kern 0.7pt 2}\,\times\\
\times\,\Bigl(1+
\frac{3\,\sqrt{2}-4}{B^{\kern 1pt 2}\,\tilde q^{\kern 0.7pt 4}}
+\frac{\tilde p\,(\sqrt{2}+1)}{2\,B^{\kern 1pt 3}
\,\tilde q^{\kern 0.7pt 7}}\Bigr)\geqslant(\sqrt{2}-1)
\,\tilde q^{\kern 0.7pt 2}\,\Bigl(1-\frac{3\,\sqrt{2}-4}
{3600^{\kern 0.7pt 4}}\,-\\
-\,\frac{\tilde p\,(\sqrt{2}+1)}{2\cdot 3600^{\kern 0.7pt 7}}\Bigr)
\approx 0.41\cdot \tilde q^{\kern 0.7pt 2}
\geqslant 0.41\cdot 3600^{\kern 0.7pt 2}\approx 5.31
\cdot 10^{\kern 0.5pt 6}>0.
\endgathered
\mytag{6.10}
$$
The left hand side of \mythetag{5.20} is treated similarly:
$$
\hskip -2em
\gathered
(\sqrt{2}-1)\,\tilde q^{\kern 0.7pt 2}
+\frac{10-7\,\sqrt{2}\,}{B^{\kern 1pt 2}\,\tilde q^{\kern 0.7pt 2}}
=(\sqrt{2}-1)\,\tilde q^{\kern 0.7pt 2}\,\Bigl(1+
\frac{3\,\sqrt{2}-4}{B^{\kern 1pt 2}\,\tilde q^{\kern 0.7pt 4}}\Bigr)\,
\geqslant\\
\geqslant(\sqrt{2}-1)
\,\tilde q^{\kern 0.7pt 2}\,\Bigl(1-\frac{3\,\sqrt{2}-4}
{3600^{\kern 0.7pt 4}}\Bigr)\geqslant 0.41\cdot 3600^{\kern 0.7pt 2}
\approx 5.31\cdot 10^{\kern 0.5pt 6}>0.
\endgathered
\mytag{6.11}
$$\par
     The root $t_4$ is delimited by the inequalities \mythetag{5.9}
and \mythetag{5.10}. Let's compare the left hand side of \mythetag{5.9}
with the right hand side of \mythetag{5.19}. Here is their difference:
$$
\Bigl((\sqrt{2}+1)\,\tilde q^{\kern 0.7pt 2}
-\frac{10+7\,\sqrt{2}\,}{B^{\kern 1pt 2}\,\tilde q^{\kern 0.7pt 2}}
\Bigr)-\Bigl((\sqrt{2}-1)\,\tilde q^{\kern 0.7pt 2}
+\frac{10-7\,\sqrt{2}\,}{B^{\kern 1pt 2}\,\tilde q^{\kern 0.7pt 2}}
\Bigr)=2\,\tilde q^{\kern 0.7pt 2}
-\frac{20}{B^{\kern 1pt 2}\,\tilde q^{\kern 0.7pt 2}}.
$$
Applying \mythetag{6.2} and \mythetag{6.3} to the above expression, 
we get the following estimate:
$$
\hskip -2em
\gathered
2\,\tilde q^{\kern 0.7pt 2}
-\frac{20}{B^{\kern 1pt 2}\,\tilde q^{\kern 0.7pt 2}}
=2\,\tilde q^{\kern 0.7pt 2}\,\Bigl(1
-\frac{10}{B^{\kern 1pt 2}\,\tilde q^{\kern 0.7pt 4}}\Bigr)\geqslant
2\,\tilde q^{\kern 0.7pt 2}\,\Bigl(1-\frac{10}{3600^{\kern 0.7pt 4}}
\Bigr)\,\approx\\
\vspace{2ex}
\approx\,1.99\cdot\tilde q^{\kern 0.7pt 2}\geqslant
1.99\cdot3600^{\kern 0.7pt 2}\approx 2.59\cdot 10^{\kern 0.5pt 7}>0.
\endgathered
\mytag{6.12}
$$
The difference of the left hand side of \mythetag{5.10} and
the right hand side of \mythetag{5.20} is treated similarly.
For this difference we have 
$$
\hskip -2em
\gathered
\Bigl((\sqrt{2}+1)\,\tilde q^{\kern 0.7pt 2}
-\frac{10+7\,\sqrt{2}\,}{B^{\kern 1pt 2}\,\tilde q^{\kern 0.7pt 2}}
-\frac{80\,\tilde p}{B^{\kern 1pt 3}\,\tilde q^{\kern 0.7pt 5}}
\Bigr)-\Bigl((\sqrt{2}-1)\,\tilde q^{\kern 0.7pt 2}\,+\\
+\,\frac{10-7\,\sqrt{2}\,}{B^{\kern 1pt 2}\,\tilde q^{\kern 0.7pt 2}}
+\frac{\tilde p}{2\,B^{\kern 1pt 3}\,\tilde q^{\kern 0.7pt 5}}
\Bigr)=
2\,\tilde q^{\kern 0.7pt 2}
-\frac{20}{B^{\kern 1pt 2}\,\tilde q^{\kern 0.7pt 2}}
-\frac{81\,\tilde p}{B^{\kern 1pt 3}\,\tilde q^{\kern 0.7pt 5}}
=\\
=2\,\tilde q^{\kern 0.7pt 2}\Bigl(1
-\frac{10}{B^{\kern 1pt 2}\,\tilde q^{\kern 0.7pt 4}}
-\frac{81\,\tilde p}{2\,B^{\kern 1pt 3}\,\tilde q^{\kern 0.7pt 7}}
\Bigr)\geqslant 2\,\tilde q^{\kern 0.7pt 2}
\Bigl(1-\frac{10}{3600^{\kern 0.7pt 4}}\,-\\
-\frac{81}{2\cdot 3600^{\kern 0.7pt 7}}\Bigr)
\approx 1.99\cdot\tilde q^{\kern 0.7pt 2}\geqslant
1.99\cdot3600^{\kern 0.7pt 2}\approx 2.59\cdot 10^{\kern 0.5pt 7}>0.
\endgathered
\mytag{6.13}
$$
The inequalities \mythetag{6.10}, \mythetag{6.11}, \mythetag{6.12}, 
and \mythetag{6.13} mean that if the inequality \mythetag{6.1} is 
fulfilled, then the asymptotic sites for the roots $t_4$ and $t_5$ 
do not intersect with each other and are located within the positive 
half-line of the imaginary axis. Summarizing this result with the above
result for the real roots $t_1$, $t_2$, $t_3$ we can formulate the 
following two theorems.\par 
\mytheorem{6.1} If\/ $\tilde p<0$ and $\tilde q\geqslant 3600\,\root{3}
\of{|\tilde p\kern 0.7pt|}$, then the equation \mythetag{2.4} has
five simple roots obeying the condition \mythetag{2.5}. Three of them
are real and positive. These three positive real roots are located 
within three disjoint asymptotic sites given by the inequalities
\mythetag{4.11}, \mythetag{4.18}, and \mythetag{4.28} respectively.
\endproclaim
\mytheorem{6.2} If\/ $\tilde p>0$ and $\tilde q\geqslant 3600\,\root{3}
\of{|\tilde p\kern 0.7pt|}$, then the equation \mythetag{2.4} has
five simple roots obeying the condition \mythetag{2.5}. Three of them
are real and positive. These three positive real roots are located 
within three disjoint asymptotic sites given by the inequalities
\mythetag{4.12}, \mythetag{4.19}, and \mythetag{4.29} respectively.
\endproclaim
\head
7. Integer points of asymptotic sites.
\endhead
     Let's consider the inequalities \mythetag{4.11} and \mythetag{4.12}
delimiting the root $t_1$ of the equation \mythetag{2.4}. Most of the 
terms in them are integer, except for following ones:
$$
\xalignat 3
&\hskip -2em
\frac{5}{B},
&&\frac{20}{B^{\kern 1pt 2}\,\tilde q^{\kern 0.7pt 2}},
&&\frac{10\,\tilde p}{B^{\kern 1pt 2}\,\tilde q^{\kern 0.7pt 3}}.
\mytag{7.1}
\endxalignat 
$$
The first term \mythetag{7.1} is optionally non-integer. If the inequalities 
\mythetag{6.2} and \mythetag{6.3} are fulfilled, then the other two terms
\mythetag{7.1} are certainly non-integer.\par
     Let's begin with the case where $5/B$ is not integer. In this case due
to \mythetag{6.3} it is separated from the nearest integer number by a distance 
not less than 1/9:
$$
\hskip -2em
\Bigl|\frac{5}{B}-n\Bigr|\geqslant \frac{1}{9}.
\mytag{7.2}
$$
The other two terms in \mythetag{7.1} are substantially smaller.
From \mythetag{6.2} and \mythetag{6.3} we get
$$
\xalignat 2
&\hskip -2em
\Bigl|\frac{20}{B^{\kern 1pt 2}\,\tilde q^{\kern 0.7pt 2}}\Bigr|
<\frac{20}{3600^{\kern 0.7pt 2}},
&&\frac{10\,\tilde p}{B^{\kern 1pt 2}\,\tilde q^{\kern 0.7pt 3}}
<\frac{10}{3600^{\kern 0.7pt 3}}.
\mytag{7.3}
\endxalignat 
$$
Combining \mythetag{7.2} and \mythetag{7.3}, we obtain the following inequality:
$$
\hskip -2em
\Bigl|\frac{5}{B}\pm \frac{20}{B^{\kern 1pt 2}\,\tilde q^{\kern 0.7pt 2}}
\pm\frac{10\,\tilde p}{B^{\kern 1pt 2}\,\tilde q^{\kern 0.7pt 3}}-n\Bigr|
\geqslant \frac{1}{9}-\frac{20}{3600^{\kern 0.7pt 2}}
-\frac{10}{3600^{\kern 0.7pt 3}}>\frac{1}{10}.
\mytag{7.4}
$$\par
     In the case where $5/B$ is integer, i\.\,e\. where $B=5$, we treat
the inequalities \mythetag{4.11} and \mythetag{4.12} separately. If $
\tilde p<0$, from \mythetag{6.2} and \mythetag{6.3} we derive
$$
\hskip -2em
-1<-\frac{20}{3600^{\kern 0.7pt 2}}-\frac{10}{3600^{\kern 0.7pt 3}}
\leqslant-\frac{20}{B^{\kern 1pt 2}\,\tilde q^{\kern 0.7pt 2}}
+\frac{10\,\tilde p}{B^{\kern 1pt 2}\,\tilde q^{\kern 0.7pt 3}}
<-\frac{20}{B^{\kern 1pt 2}\,\tilde q^{\kern 0.7pt 2}}<0.
\mytag{7.5}
$$
If $\tilde p>0$ we need an additional condition for $\tilde p$ and $\tilde q$:
$$
\hskip -2em
2\,\tilde q>|\tilde p\kern 0.7pt|
\mytag{7.6}
$$
Provided the condition \mythetag{7.6} is fulfilled, from 
\mythetag{6.2} and \mythetag{6.3} we derive
$$
\hskip -2em
-1<-\frac{20}{3600^{\kern 0.7pt 2}}\leqslant
-\frac{20}{B^{\kern 1pt 2}\,\tilde q^{\kern 0.7pt 2}}
<-\frac{20}{B^{\kern 1pt 2}\,\tilde q^{\kern 0.7pt 2}}
+\frac{10\,\tilde p}{B^{\kern 1pt 2}\,\tilde q^{\kern 0.7pt 3}}<0.
\mytag{7.7}
$$
If the condition \mythetag{7.6} is not fulfilled, then from 
\mythetag{6.2} and \mythetag{6.3} we derive 
$$
-1<-\frac{20}{3600^{\kern 0.7pt 2}}\leqslant
-\frac{20}{B^{\kern 1pt 2}\,\tilde q^{\kern 0.7pt 2}}
<-\frac{20}{B^{\kern 1pt 2}\,\tilde q^{\kern 0.7pt 2}}
+\frac{10\,\tilde p}{B^{\kern 1pt 2}\,\tilde q^{\kern 0.7pt 3}}
<\frac{20}{3600^{\kern 0.7pt 2}}
+\frac{10}{3600^{\kern 0.7pt 3}}<1.
\quad
\mytag{7.8}
$$
\mytheorem{7.1}If $\tilde p\neq 0$, then for each $\tilde q\geqslant 3600
\,\root{3}\of{|\tilde p\kern 0.7pt|}$ and for each $B\neq 5$ in 
\mythetag{4.6} the asymptotic site given by the inequalities \mythetag{4.11} 
and \mythetag{4.12} has no integer points.
\endproclaim
Theorem~\mythetheorem{7.1} is proved by applying \mythetag{7.4} to
\mythetag{4.11} and \mythetag{4.12}.\par
\mytheorem{7.2}If $\tilde p<0$ and $B=5$, then for each $\tilde q\geqslant 
3600\,\root{3}\of{|\tilde p\kern 0.7pt|}$ the asymptotic site given by the 
inequalities \mythetag{4.11} has no integer points.
\endproclaim
Theorem~\mythetheorem{7.2} is proved by applying \mythetag{7.5} to
\mythetag{4.11}.\par
\mytheorem{7.3}If $\tilde p>0$ and $B=5$, then for each $\tilde q$ 
such that $\tilde q\geqslant 3600\,\root{3}\of{|\tilde p\kern 0.7pt|}$ and 
$2\,\tilde q>\tilde p$ the asymptotic site given by the 
inequalities \mythetag{4.12} has no integer points.
\endproclaim
Theorem~\mythetheorem{7.3} is proved by applying \mythetag{7.7} to
\mythetag{4.12}.\par
\mytheorem{7.4}If $\tilde p\neq 0$ and $B=5$, then for $\tilde q\geqslant 3600\,\root{3}\of{|\tilde p\kern 0.7pt|}$ the asymptotic site given by the 
inequalities \mythetag{4.11} and \mythetag{4.12} has at most one integer 
point given by the formula $t=25\,\tilde q^6+10\,\tilde q^4-10\,\tilde p
\,\tilde q^3-2\,\tilde q^2-2\,\tilde p\,\tilde q+\tilde p^2+1$, where
$\tilde p>0$. 
\endproclaim
Theorem~\mythetheorem{7.4} is proved by applying \mythetag{7.8} to
\mythetag{4.12}.\par
     Note that the inequalities \mythetag{4.18} and \mythetag{4.19}
are quite similar to \mythetag{4.11} and \mythetag{4.12}. Applying
the estimates \mythetag{7.4}, \mythetag{7.5}, \mythetag{7.7}, and 
\mythetag{7.8} to them we can prove the following four theorems for
the corresponding asymptotic site.\par
\mytheorem{7.5}If $\tilde p\neq 0$, then for each $\tilde q\geqslant 3600
\,\root{3}\of{|\tilde p\kern 0.7pt|}$ and for each $B\neq 5$ in 
\mythetag{4.6} the asymptotic site given by the inequalities \mythetag{4.18} 
and \mythetag{4.19} has no integer points.
\endproclaim
\mytheorem{7.6}If $\tilde p>0$ and $B=5$, then for each $\tilde q\geqslant 
3600\,\root{3}\of{|\tilde p\kern 0.7pt|}$ the asymptotic site given by the 
inequalities \mythetag{4.19} has no integer points.
\endproclaim
\mytheorem{7.7}If $\tilde p<0$ and $B=5$, then for each $\tilde q$ 
such that $\tilde q\geqslant 3600\,\root{3}\of{|\tilde p\kern 0.7pt|}$ and 
$2\,\tilde q>|\tilde p\kern 0.7pt|$ the asymptotic site given by the 
inequalities \mythetag{4.18} has no integer points.
\endproclaim
\mytheorem{7.8}If $\tilde p\neq 0$ and $B=5$, then for $\tilde q\geqslant 3600\,\root{3}\of{|\tilde p\kern 0.7pt|}$ the asymptotic site given by the 
inequalities \mythetag{4.18} and \mythetag{4.19} has at most one integer 
point given by the formula $t=25\,\tilde q^6-10\,\tilde q^4-10\,\tilde p
\,\tilde q^3-2\,\tilde q^2+2\,\tilde p\,\tilde q+\tilde p^2-1$, where
$\tilde p<0$. 
\endproclaim
     Let's proceed to the inequalities \mythetag{4.28} and \mythetag{4.29}
and let's begin with the case where $16/B$ is not integer. In this case due
to \mythetag{6.3} the number $16/B$ is separated from the nearest integer 
number by a distance not less than 1/9:
$$
\hskip -2em
\Bigl|\frac{16}{B}-n\Bigr|\geqslant \frac{1}{9}.
\mytag{7.9}
$$
From \mythetag{6.2} and \mythetag{6.3} for the other fractional term in 
\mythetag{4.28} and \mythetag{4.29} we derive 
$$
\hskip -2em
\Bigl|\frac{32\,\tilde p}{B^{\kern 1pt 2}\,
\tilde q^{\kern 0.7pt 3}}\Bigr|\leqslant\frac{32}{3600^{\kern 0.7pt 3}}
\mytag{7.10}
$$
Combining \mythetag{7.9} and \mythetag{7.10}, we obtain the following 
inequality:
$$
\hskip -2em
\Bigl|\frac{16}{B}\pm \frac{32\,\tilde p}{B^{\kern 1pt 2}\,
\tilde q^{\kern 0.7pt 3}}-n\Bigr|\geqslant \frac{1}{9}
-\frac{32}{3600^{\kern 0.7pt 3}}\geqslant\frac{1}{10}.
\mytag{7.11}
$$
The case where $16/B$ is integer is more simple. In this case we write
the inequality \mythetag{7.10} in the following slightly modified form:
$$
\hskip -2em
0<\Bigl|\frac{32\,\tilde p}{B^{\kern 1pt 2}\,
\tilde q^{\kern 0.7pt 3}}\Bigr|\leqslant
\frac{32}{3600^{\kern 0.7pt 3}}<1.
\mytag{7.12}
$$
In both cases, applying either \mythetag{7.11} or \mythetag{7.12}
to the inequalities \mythetag{4.28} and \mythetag{4.29}, we can 
prove the following theorem. 
\mytheorem{7.9}If $\tilde p\neq 0$, then for each $\tilde q\geqslant 3600
\,\root{3}\of{|\tilde p\kern 0.7pt|}$ and for each $B$ in \mythetag{4.6} 
the asymptotic site given by the inequalities \mythetag{4.28} 
and \mythetag{4.29} has no integer points.
\endproclaim
\head
The exceptional case. 
\endhead
     As we noted above the case $\tilde p=0$ is exceptional, see
Theorems~\mythetheorem{7.1} through \mythetheorem{7.9}. This case should be
studied separately. Let's substitute $\tilde p=0$ into the equation
\mythetag{2.4}. The resulting equation can be written explicitly:
$$
\pagebreak 
\gathered
t^{10}+(6\,\tilde q^{\kern 0.7pt 4}-2\,B^{\kern 1pt 4}
\,\tilde q^{\kern 0.7pt 12}-\tilde q^{\kern 0.7pt 8}\,B^{\kern 1pt 2})\,t^8
+(B^{\kern 1pt 8}\,\tilde q^{\kern 0.7pt 24}+10\,\tilde q^{\kern 0.7pt 12}
\,B^{\kern 1pt 2}+4\,\tilde q^{\kern 0.7pt 16}\,B^{\kern 1pt 4}\,-\\
-\,14\,\tilde q^{\kern 0.7pt 20}\,B^{\kern 1pt 6}
+\tilde q^{\kern 0.7pt 8})\,t^6+(14\,\tilde q^{\kern 0.7pt 20}
\,B^{\kern 1pt 4}-4\,\tilde q^{\kern 0.7pt 24}\,B^{\kern 1pt 6}
-\tilde q^{\kern 0.7pt 16}\,B^{\kern 1pt 2}
-\tilde q^{\kern 0.7pt 32}\,B^{\kern 1pt 10}\,-\\
-\,10\,\tilde q^{\kern 0.7pt 28}\,B^{\kern 1pt 8})\,t^4
+(2\,\tilde q^{\kern 0.7pt 28}\,B^{\kern 1pt 6}-6\,\tilde q^{\kern 0.7pt 36}
\,B^{\kern 1pt 10}+\tilde q^{\kern 0.7pt 32}\,B^{\kern 1pt 8})\,t^2
-\tilde q^{\kern 0.7pt 40}\,B^{\kern 1pt 10}=0.
\endgathered\quad
\mytag{8.1}
$$
\mytheorem{8.1} For each $B$ in \mythetag{4.6} the polynomial in the
left hand side of the equation \mythetag{8.1} is irreducible in the ring 
$\Bbb Z[t]$. 
\endproclaim
Theorem~\mythetheorem{8.1} is proved by means of direct computations. It
means that for each $B$ in \mythetag{4.6} the equation \mythetag{8.1}
has no integer roots. 
\head
9. Application to the cuboid problem.
\endhead
     Theorems~\mythetheorem{7.1} through \mythetheorem{7.9} are based on
the inequality \mythetag{6.1}, where $\tilde p\neq 0$. Theorem~\mythetheorem{8.1}
means that we can omit the condition $\tilde p\neq 0$. Transforming 
\mythetag{6.1} back to the initial variables $p$ and $q$ with the use of 
\mythetag{2.1}, we get the inequalities 
$$
\hskip -2em
B\,q^{\kern 0.7pt 3}-\frac{q^{\kern 0.7pt 3}}{3600^{\kern 0.7pt 3}}
\leqslant p\leqslant B\,q^{\kern 0.7pt 3}
+\frac{q^{\kern 0.7pt 3}}{3600^{\kern 0.7pt 3}}.
\mytag{9.1}
$$
From Theorems~\mythetheorem{6.2}, \mythetheorem{7.1}, \mythetheorem{7.5},
and \mythetheorem{7.9} we derive the following result. 
\mytheorem{9.1} For each $B\neq 5$ in \mythetag{4.6} if the 
inequalities \mythetag{9.1} are fulfilled, then the tenth degree cuboid 
characteristic equation \mythetag{1.1} produces no perfect cuboids. 
\endproclaim
     The case $B=5$ is special. In this case we have the additional 
condition \mythetag{7.6}. Upon transforming \mythetag{7.6} back to the 
initial variables $p$ and $q$ it looks like
$$
\hskip -2em
B\,q^{\kern 0.7pt 3}-2\,q<p< B\,q^{\kern 0.7pt 3}+2\,q.
\mytag{9.2}
$$
From Theorems~\mythetheorem{6.2}, \mythetheorem{7.2}, \mythetheorem{7.3}, 
\mythetheorem{7.6}, \mythetheorem{7.7}, and \mythetheorem{7.9} we derive 
the following result. 
\mytheorem{9.2} For $B=5$ if the inequalities \mythetag{9.1} and 
\mythetag{9.2} are fulfilled, then the tenth degree cuboid 
characteristic equation \mythetag{1.1} produces no perfect cuboids. 
\endproclaim
     The inequalities \mythetag{9.2} do not follow from \mythetag{9.1}.
They become very restrictive for large $q$ as $q\to+\infty$. 
Theorems~\mythetheorem{7.4} and \mythetheorem{7.8} can be applied in order
to remove this restriction. However, they do not change the state of 
affairs in general. Therefore this step is left for one of the future 
papers. 
\head
10. Conclusions.
\endhead
    The main result of this paper is presented by Theorems~\mythetheorem{9.1}
and \mythetheorem{9.2}. Theorems~\mythetheorem{9.1} and \mythetheorem{9.2} 
shrink the nonlinear region on the $p\,q$\,-\,coordinate plane by cutting off 
nine narrow strips surrounding nine cubic parabolas
$$
p=B\,q^{\kern 0.7pt 3}\text{, \ where \ }B=1,\,2,\,\ldots,\,9.
$$
These strips are outlined by the inequalities \mythetag{9.1} and \mythetag{9.2}.
They are annexed to the no cuboid region (see Fig\.~1.1) thus reducing the area 
where perfect cuboids are still potentially possible\footnotemark.\par  
\footnotetext{\ The paper \mycite{70} has been recently published in ArXiv. It 
says that there are no perfect cuboids. However this paper is not yet verified
by the mathematical community. Therefore alternative approaches to the perfect
cuboid problem can be developed for some while.}
\head
11. Acknowledgments
\endhead
     On May 19, 2015, I have reported the papers \myciterange{1}{1}{--}{3} in  
the Ufa all-city seminar on differential equations of mathematical physics named 
after A\.~M\. Ilyin. \pagebreak This seminar brings together many experts in the 
field of asymptotics residing in our city. I am grateful to L\.~A\.~Kalyakin, the 
chairman of the seminar, for the opportunity to give my talk. I am also grateful 
to V\.~Yu\.~Novokshenov, the other chairman, and to all participants  of the 
seminar for their attention and comments.\par 
\adjustfootnotemark{-1}

\enddocument
\end